\documentclass{article}   	
\usepackage{geometry}               
\usepackage{graphicx}
\usepackage[usenames,dvipsnames]{xcolor}
\usepackage{amssymb,amsmath,amsthm,amsfonts}

\allowdisplaybreaks

\usepackage[english]{babel}
\usepackage[utf8]{inputenc}
\usepackage[colorlinks]{hyperref}
\hypersetup{linkcolor=blue,citecolor=blue,filecolor=black,urlcolor=blue}

\usepackage{mathtools}
\usepackage{dsfont}

\usepackage[sc]{mathpazo}

\newtheorem{proposition}{Proposition}[section]
\newtheorem{theorem}[proposition]{Theorem}

\newtheorem{lemma}[proposition]{Lemma}

\theoremstyle{definition}

\theoremstyle{remark}
\newtheorem{remark}[proposition]{Remark}
\numberwithin{equation}{section}

\newcommand{\eps}{\varepsilon}

\newcommand{\N}{{\mathbb{N}}}

\newcommand{\R}{{\mathbb{R}}}

\newcommand{\bm}{\mathbf}
\newcommand{\bx}{{\bm{x}}}
\newcommand{\by}{{\bm{y}}}

\newcommand{\Lcal}{{\mathcal{L}}}
\newcommand{\Mcal}{{\mathcal{M}}}

\DeclareMathOperator{\diam}{diam}
\DeclareMathOperator{\supp}{supp}

\newcommand{\ind}[1]{\chi_{#1}}

\title{Asymptotic properties of an optimal principal eigenvalue with spherical weight and Dirichlet boundary conditions}
\author{Lorenzo Ferreri, Gianmaria Verzini}

\begin{document}
\maketitle

\begin{abstract}
We consider a weighted eigenvalue problem for the Dirichlet laplacian in a 
smooth bounded domain $\Omega\subset \R^N$, where the bang-bang weight 
equals a positive constant $\overline{m}$ on a ball $B\subset\Omega$ and a 
negative constant $-\underline{m}$ on $\Omega\setminus B$. The 
corresponding positive principal eigenvalue provides a threshold to detect 
persistence/extinction of a species whose evolution is described by 
the heterogeneous Fisher-KPP equation in population dynamics. In 
particular, we study the minimization of such eigenvalue with respect to the position of $B$ in $\Omega$. We provide sharp asymptotic expansions 
of the optimal eigenpair in the singularly perturbed regime in which the volume 
of $B$ vanishes. We deduce that, up to subsequences, the optimal ball 
concentrates at a point maximizing the distance from $\partial\Omega$.
\end{abstract}
\noindent
{\footnotesize \textbf{AMS-Subject Classification}}. 
{\footnotesize 49R05, 92D25, 47A75, 35B40
}

\section{Introduction}

Let $\Omega \subset \R^{N}$ denote an open, bounded and connected set with smooth 
boundary (for instance, $\partial\Omega$ of class $C^2$ is enough), and  
$m\in L^\infty(\Omega)$ be a sign-changing weight having non-trivial positive part. 
Throughout this paper we denote with $\lambda^1(m,\Omega)$ the positive principal eigenvalue of the Dirichlet problem
\begin{equation}\label{GeneralDiffProblem}
\begin{cases}
-\Delta u = \lambda m u & \text{in } \Omega \, , \\
u = 0 & \text{on } \partial\Omega 
\end{cases}
\end{equation}
(actually, we will simply write $\lambda^1(m)$ unless $\Omega$ is replaced by some 
other domain). Precisely,
\[
\lambda^1(m) = \inf_{\substack{u \in H^1_0(\Omega) \\ \int_{\Omega}m u^2 > 0}} 
\frac{\int_{\Omega}|\nabla u|^2}{\int_{\Omega}m u^2}.
\]
Moreover, it is well known that $\lambda^1(m)$ is simple, it is a continuous function 
of $m$ in the $L^{N/2}(\Omega)$-norm, and it is achieved by an eigenfunction 
$u \in H^1_0(\Omega)$ that can be chosen strictly positive over $\Omega$ (see e.g.~\cite{deFiguereido:EqDiff82}).

The analysis of the dependence of $\lambda^1$ on $m$ is relevant for different reasons, 
both theoretical and related to applications. In particular, here we are interested in 
its role on the persistence of a species in population dynamics, as first enlighten 
by Cantrell and Cosner in some seminal papers 
\cite{CantrellCosner:indefiniteWeight,MR1112065,MR1105497}. Indeed, let us consider a 
logistic reaction-diffusion model for the dispersal of a species in an heterogeneous 
environment, described by the weight $m$ (which is positive/negative in favourable/
hostile zones of the habitat, respectively). In such setting, Dirichlet 
conditions entail an environment surrounded by uninhabitable regions. Then it is well 
known that the species persists, i.e.\ it evolves towards a non-trivial steady state, 
if and only if $\lambda^1(m)$ is below a certain explicit threshold, 
see~\cite{MR2191264,Berestycki:PeriodicallyFragmented}. As a consequence, different 
optimization problems arise concerning $\lambda^1$, see 
\cite{ly,RoquesHamel:OptimalHabitat,dego,Lamboley:OptimizersRobin,MR3771424,mapeve,dipierro2021nonlocal,Verzini:Neumann} and references therein, also for different boundary 
conditions and/or operators.

As the optimization problem we deal with in this paper is related to that originally 
considered in \cite{CantrellCosner:indefiniteWeight}, let us first describe the 
latter. The problem considered by Cantrell and Cosner is to minimize $\lambda^1(m)$, 
with $m$ belonging to the admissible class 
\[
\mathcal{M}_{m_0} \coloneqq \{ m \in L^{\infty}(\Omega) : -\underline{m} \le m \le \overline{m},\ \int_{\Omega}m = m_0\Lcal(\Omega),\ \mathcal{L}(\{m>0\}) >0 \}
\] 
for some fixed parameters $\underline{m},\overline{m},m_0$ with 
$-\underline{m} < m_0 < \overline{m}$ (we denote with $\Lcal(A)$ the Lebesgue measure 
of the measurable set $A$). It is shown in 
\cite[Thm. 3.9]{CantrellCosner:indefiniteWeight} that the optimization problem
$\inf_{\Mcal}\lambda^1(m)$ is achieved by a so called bang-bang weight, i.e. by a 
weight belonging to the set
\[
\mathcal{BB}_{\varepsilon} \coloneqq  \{ \overline{m} \ind{E}-\underline{m}\ind{\Omega\setminus E}: \Lcal(E)=\eps \},
\qquad\text{with}\qquad 
0<\eps\coloneqq \frac{m_0+\underline{m}}{\overline{m}+\underline{m}}\Lcal(\Omega)
<\Lcal(\Omega),
\]
therefore the true unknown in the above minimization problem is the shape of the 
optimal favorable region $E$. In other words, with a little abuse of notation we 
write
\[
\lambda^1(\overline{m} \ind{E}-\underline{m}\ind{\Omega\setminus E}) = 
\lambda^1(E,\Omega) =
\lambda^1(E)
\]
(again, we drop the dependence on the box $\Omega$ when no confusion arises) and we obtain that the minimization problem considered by Cantrell and Cosner is
\begin{equation}\label{GeneralProblem}
\min_{m \in \mathcal{M}_{m_0}} \lambda^1(m) = \min_{m \in \mathcal{BB}_{\varepsilon}} \lambda^1(m)= \min_{\Lcal(E)=\eps} \lambda^1(E).
\end{equation}

About the optimal set $E$, the proof of the above equivalence rests on the 
fact that it is a superlevel set of the associated positive eigenfunction. Apart from 
this, not much information is available, and natural questions concern its qualitative 
properties (connectedness, symmetry) and notably its location inside $\Omega$. A complete answer is available 
only when $\Omega\subset \R^N$, $N\ge1$, is a ball: in this case, by symmetrization 
arguments, $E\subset\Omega$ is a concentric ball \cite[Rmk. 3.10]{CantrellCosner:indefiniteWeight}. We aim at extracting further information about the general problem, in the 
singularly perturbed regime $\eps\to0$. In this case, based on the analysis performed 
in \cite{Verzini:Neumann} for the analogous problem with Neumann boundary conditions, 
we expect that the optimal set is connected and that it concentrates somewhere in 
$\Omega$, with a shape which is asymptotically spherical. On the other hand, the 
techniques in \cite{Verzini:Neumann} seem not enough to localize the point of 
concentration. For this reason we deal here with a reduced version of 
\eqref{GeneralProblem}, characterized by the additional assumption that $E$ is a 
sphere. Our main goals concern the analysis of such reduced problem, as $\eps\to0^+$, providing 
sharp asymptotics of the optimal eigenvalues/eigenfunctions, and in particular the 
precise position of the points in which the optimal sets concentrate.

To describe the problem we consider, let us write
\[
B^{\varepsilon}(\mathbf{x}) = B_{r(\varepsilon)}(\mathbf{x}), \qquad \text{where 
$r(\varepsilon)$ is such that }
\Lcal(B^\eps) = \eps.
\]  
Moreover, we introduce the compact set 
\[
\hat{\Omega}_{\varepsilon} \coloneqq \{ \mathbf{x} \in \Omega : d(\mathbf{x}, \partial\Omega) \ge  r(\varepsilon) \}.
\]
Of course this set is non-empty, for $0<\eps\le\eps_0$ sufficiently small, 
and $\mathbf{x}\in\hat\Omega$ implies $B^{\varepsilon}(\mathbf{x}) \subset \Omega$. 

We are interested in the study of
\begin{equation}\label{sphericalProblem}
\lambda_{\varepsilon} \coloneqq  \inf_{\mathbf{x} \in \hat{\Omega}_{\varepsilon} } \lambda^1(B^{\varepsilon}(\mathbf{x})) \qquad \text{if } \varepsilon \to 0^{+}.
\end{equation}
Since the map $\bx\mapsto \lambda^1(B^{\varepsilon}(\mathbf{x}))$ is continuous in $L^p$, it is immediate to prove the following.
\begin{proposition}\label{prop:PropositionMinimum}
For any $\varepsilon > 0$ there exixts $\mathbf{x}_{\varepsilon} \in \hat\Omega_\eps$ such that
\[
\lambda_{\varepsilon} = \lambda^1(B^{\varepsilon}(\mathbf{x}_{\varepsilon})) ,
\]
achieved by an eigenfunction $u_{\varepsilon} \in H^1_0(\Omega)$ which we assume 
positive and normalized in $L^2(\Omega)$ (unless otherwise stated).
\end{proposition}

We remark that, contrarily to \eqref{GeneralProblem}, problem \eqref{sphericalProblem}, in general, lacks the strong property that $B^{\varepsilon}(\mathbf{x}_{\varepsilon})$ is a superlevel set for $u_{\varepsilon}$. However, this lack of information is compensated by the fact that, by the maximum principle, the maximum points of 
$u_\eps$ belong to $B^{\varepsilon}(\mathbf{x}_{\varepsilon})$.

We study the properties of the family $\{ u_{\varepsilon} \}$ for $\varepsilon \to 0^{+}$ through a blow-up argument. More precisely, we introduce the quantities
\begin{equation}\label{eq:kbeta_intro}
k_{\varepsilon} \coloneqq \varepsilon^{1/N},\qquad \beta_\eps = \frac1{k_{\varepsilon}} \coloneqq \varepsilon^{-1/N},
\end{equation}
and the blow-up sequences
\begin{equation}\label{eq:blowupseq_intro}
\tilde{u}_\eps(\mathbf{x}) \coloneqq k_\eps^{N/2}  u_{\varepsilon}(\mathbf{x}_{\varepsilon} + k_\eps \mathbf{x}) \qquad \text{in}\qquad \tilde{\Omega}_\eps \coloneqq \left\{ \mathbf{x} \in \R^N : \mathbf{x}_{\varepsilon} + k_\eps \mathbf{x} \in \Omega \right\},
\end{equation}
so that $\|\tilde u_\eps\|_{L^2(\tilde \Omega_\eps)} =1$. 
The key point in this analysis is that, by rescaling, all the optimal 
weights coincide. Whence, the functions 
$\tilde{u}_\eps$ solve, in $H_0^1(\tilde{\Omega}_\eps)$:
\begin{equation}\label{GeneralDiffProblemBU_intro}
\begin{cases}
-\Delta \tilde{u}_\eps = \tilde{\lambda}_{\eps} \tilde{m}_0 \tilde{u}_\eps & \text{in } \tilde{\Omega}_\eps \, , \\
\tilde{u}_\eps = 0 & \text{on } \partial\tilde{\Omega}_\eps \, ,
\end{cases}
\end{equation}
with 
\begin{equation}\label{eq:lambdabu_intro}	
\tilde{\lambda}_{\eps} \coloneqq k_\eps^2\lambda_{\varepsilon}
\qquad\text{and}\qquad
\tilde{m}_0 = \overline{m} \ind{B^1(\mathbf{0})}-\underline{m}\ind{ \R^N \setminus B^1(\mathbf{0}) }. 
\end{equation}
We are naturally led to consider the corresponding limit problem on $\R^N$, which was
already considered in \cite{Verzini:Neumann}. In particular, we denote with $w \in H^1(\R^N)$ the only positive, radial and $L^2(\R^N)$ normalized solution of
\begin{equation}\label{problemOverRn}
-\Delta w = \tilde{\lambda}_0 \tilde{m}_0 w \quad \text{in } \R^N \, ,
\end{equation}
where $\tilde{\lambda}_0$ is the positive principal eigenvalue of 
\eqref{problemOverRn}, which exists finite. In particular, $w$ is radially 
decreasing and vanishes exponentially at infinity. 
Then, we prove the following.
\begin{proposition}\label{prop:PropositionBlowUp}
For any vanishing sequence $\varepsilon_{n}$ there exists a subsequence still denoted with $\varepsilon_{n}$ such that: 
\begin{itemize}
\item[(i)] $\lambda_{\varepsilon_{n}} \to \lambda_0$,
\item[(ii)] $\tilde u_{\varepsilon_{n}} \to w$ in $C^{1, \alpha}(K)$ for $n \to +\infty$, for any $0 < \alpha < 1$ and $K\subset\R^N$ compact set,
\item[(iii)] $\tilde u_{\varepsilon_{n}} \to w$ in $H^1(\R^N)$ for $n \to +\infty$.
\end{itemize}
\end{proposition}

Going back to the original variables, the above proposition yields, up to subsequences,  
concentration of the optimal eigenfunctions at some point which is a limit of the 
centers of the optimal favorable sets. To proceed with our analysis, we need to sharpen 
the above asymptotics to obtain a 
further term in the expansion. This will lead to the precise localization of the 
concentration 
point. 

Our strategy is inspired by the techniques presented in 
\cite{NiWei:Spike95}, where Ni and Wei analyze the concentration of the least energy 
solutions for semilinear Dirichlet problems. Even though our problem is quite 
different, an important analogy with their analysis is that the limit solutions of both problems decay at infinity with the same rate. More 
precisely, let us denote with $P_{\tilde{\Omega}_{\varepsilon}}w$ the 
$H^1_0(\tilde{\Omega}_{\varepsilon})$-projection of $w$ and with 
$\tilde{\Psi}_{\varepsilon}$ the solution to
\[
\begin{cases}
- k_{\varepsilon} \Delta \tilde{\Psi}_{\varepsilon} + |\nabla \tilde{\Psi}_{\varepsilon}|^2 - \tilde{\lambda}_0 \underline{m} = 0 & \text{in } \Omega \; , \\
\tilde{\Psi}_{\varepsilon} = -k_{\varepsilon} \log(w)(\frac{\mathbf{x}-\mathbf{x}_{\varepsilon}}{k_{\varepsilon}}) & \text{on }  \partial\Omega \; ,
\end{cases}
\]
which is unique by the Comparison Principle for Quasilinear Equations, see e.g.  \cite[Thm. 10.1]{GilbargTrudinger:Elliptic98}. We prove the following.
\begin{theorem}\label{theorem:TheoremAsymptExpansion}
For any vanishing sequence $\varepsilon_{n}$ there exists a subsequence still denoted with $\varepsilon_{n}$ such that, for  $n \to +\infty$: 
\begin{itemize}
\item[(i)] $\displaystyle d(\mathbf{x}_{\varepsilon_n}, \partial\Omega) \to \max_{\mathbf{p}\in\Omega}d(\mathbf{p}, \partial\Omega)
$;
\item[(ii)]
$\displaystyle
\tilde{\Psi}_{\varepsilon_{n}}(\mathbf{x}_{\varepsilon_{n}}) \to 2 \sqrt{\tilde{\lambda}_0 \underline{m}} \, \max_{\mathbf{p}\in\Omega}d(\mathbf{p},  \partial\Omega)
$;
\item[(iii)] 
$\displaystyle
\lambda_{\varepsilon_n} = \frac{1}{k_{\varepsilon_n}^2} \left( \tilde{\lambda}_0 + \Phi \, e^{- \beta_{\varepsilon_n} \tilde{\Psi}_{\varepsilon_n}(\mathbf{x}_{\varepsilon}) } + o\left(e^{- \beta_{\varepsilon_n} \tilde{\Psi}_{\varepsilon_n}(\mathbf{x}_{\varepsilon}) }\right) \right)$,
where $\Phi$ is a positive constant.
\item[(iv)] $\displaystyle u_{\varepsilon_n} = P_{\tilde{\Omega}_{\varepsilon_n}}w + e^{- \beta_{\varepsilon_n} \tilde{\Psi}_{\varepsilon_n}(\mathbf{x}_{\varepsilon_n}) } \phi_{\varepsilon_n}$, with $\phi_{\varepsilon_n} \in H^2(\tilde\Omega_{\varepsilon_n})$ and $\| \phi_{\varepsilon_n} \|_{H^2(\tilde\Omega_{\varepsilon_n})} \le C$ uniformly for $n \to +\infty$. 
\end{itemize}
\end{theorem}

A straightforward consequence of our result is that, for our model problem 
\eqref{sphericalProblem}, the optimal set concentrates as $\eps\to0^+$, at points 
of $\Omega$ which are at maximal distance from the boundary.

The next step in our analysis will be to analyze problem \eqref{GeneralProblem} as 
$\eps\to0$. It is natural to conjecture that one should be able to extend at least 
part of our results also to this setting. Actually, the main obstruction in this 
direction is that, for the general problem, the blow-up sequences satisfy a problem 
similar to \eqref{GeneralDiffProblemBU_intro}, with a weight $\tilde m_\eps$ which 
is no more independent of $\eps$. Nonetheless, on the one hand, it should be 
possible to obtain a result similar to Proposition \ref{prop:PropositionBlowUp} by 
combining our techniques with those exploited in \cite{Verzini:Neumann} for the 
Neumann problem. In particular, this should imply that the rescaled optimal sets
$\tilde E_\eps$ converge to $B^1$ in the blow-up scale. On the other hand, the 
extension of Theorem \ref{theorem:TheoremAsymptExpansion} is much more difficult, as it 
requires information on the velocity at which the above set convergence happens. This 
analysis will be the object of further studies.

This paper is structured as follows: in Section \ref{section:minimizer} we recall some 
preliminary facts, in particular about features of problem \eqref{problemOverRn} (part 
of which will be proved in the appendix); in Sections \ref{blowUp} and 
\ref{section:H1Convergence} we carry out the first part of the blow-up analysis, 
proving Proposition \ref{prop:PropositionBlowUp}; finally, Sections \ref{section:EstimateAbove} and \ref{section:EstimateBelow} are devoted to the quantitative study of problem \eqref{sphericalProblem}, leading to the proof 
of Theorem \ref{theorem:TheoremAsymptExpansion}. 

\subsection*{Notation}

\begin{itemize}
\item $\Lcal(A)$ is the Lebesgue measure of the measurable set $A$;
\item $B_r(\mathbf{x})$ is the ball of center $\bx$ and radius $r$; $B^{\varepsilon}(\mathbf{x})$ is the ball of center $\bx$ and measure $\eps$;
\item for subsequences $k_{\eps_n},u_{\eps_n},\lambda_{\eps_n},\dots$ we simply write 
$k_{n},u_{n},\lambda_{n},\dots$
\item $C,C'$ and so on denote non-negative universal constants, which we need 
not to specify, and which may vary from line to line.
\end{itemize}

\section{Preliminaries}\label{section:minimizer}

In this section we recall some basic properties of the solution $w\in H^1(\R^N)$ 
of the limit problem \eqref{problemOverRn}. Most of them are contained in 
\cite[Sect. 2]{Verzini:Neumann}, in particular Theorem 2.2 therein (up to some rescaling). It is convenient 
to split such results in two parts. First, we address the variational 
characterization of $w$ and its connection with optimization 
problems on spheres. 

\begin{proposition}\label{result:ResultProblemOverRn}
Consider the class 
\[
\mathcal{M}' \coloneqq \{ m \in L^{\infty}(\R^N): -\underline{m} \le m \le \overline{m},\, \int_{\Omega}(m+\underline{m}) \le \underline{m} + \overline{m} \}
\]
The quantity
\begin{equation}\label{lambda0defRn}
\tilde{\lambda}_0 \coloneqq \inf_{m \in \mathcal{M}'} \inf_{\substack{u \in H^1(\R^N) \\ \int_{\R^N}m u^2 > 0}} \frac{\int_{\R^N}|\nabla u|^2}{\int_{\R^N}m u^2}  \, ,
\end{equation}
is a positive minimum, uniquely attained, up to a translation, by 
\[
\tilde{m}_0 \coloneqq \overline{m} \ind{B^1(\mathbf{0})}-\underline{m}\ind{ \R^N \setminus B^1(\mathbf{0}) }
\]
and $w \in H^1(\R^N)$ that is the only solution in $H^1(\R^N)$, up to a multiplicative constant, of
\begin{equation}\label{eqn:sulutionRn}
-\Delta w = \tilde{\lambda}_0 \tilde{m}_0 w \quad \text{in } \R^N \, .
\end{equation}
The function $w$ has constant sign and is radially symmetric with respect to the origin. In particular, if $w$ is chosen to be positive, it is radially strictly decreasing.%

Moreover, $\tilde{\lambda}_0$ is characterized by
\begin{equation}\label{eqn:characterizationTildeLambda0}
\tilde{\lambda}_0 = \lim_{n \to +\infty}\min_{B^1(\mathbf{x})\subset 
B^{A_n}(\mathbf{0})} \lambda^1(B^1(\mathbf{x}),B^{A_n}(\mathbf{0})) = \lim_{n \to +\infty} \lambda^1(B^1(\mathbf{0}),B^{A_n}(\mathbf{0})),
\end{equation}
for every diverging sequence $(A_n)_n$. 
\end{proposition}

Most of Proposition \ref{result:ResultProblemOverRn} is contained in \cite[Sect. 2]{Verzini:Neumann}. The only thing left to prove is property \eqref{eqn:characterizationTildeLambda0}: on the one hand, the inequality 
$\tilde{\lambda}_0\le  \lambda^1(B^1(\mathbf{0}),B^{A_n}(\mathbf{0}))$ holds true for every $A_n$ large since $H^1_0(B^{A_n})\subset H^1(\R^N)$ (by trivial extension); the opposite inequality easily follows by approximating in $H^1(\R^N)$ the function $w$ with a sequence $w_k$ of $C_c^{\infty}$ functions with expanding compact supports.

Since the radial symmetry of (at least) any non-negative function in $H^1(\R^N)$ attaining $\tilde{\lambda}_0$ in \eqref{lambda0defRn} is the key to prove Proposition \ref{prop:PropAwayFromBoundary}, for the sake of completeness its proof is recalled in the appendix. 

To conclude this discussion, we remark that $w$ is actually explicit, 
in terms of Bessel functions. In particular, this provides its behavior at infinity 
which is crucial both in the proof of Proposition \ref{prop:PropositionBlowUp} and 
in the extension of the techniques in \cite{NiWei:Spike95}.

\begin{proposition}\label{result:ResultShapeSolutionOverRn}
Let $w$ be the function defined in Proposition \ref{result:ResultProblemOverRn}. If expressed in radial coordinates $w = w(r)$, it is the only positive,  $L^2(\R^N)$ normalized solution of
\begin{equation}
\begin{cases}\label{eq:ODEsystem1}
w''(r) + \frac{N-1}{r}w'(r)+\tilde{\lambda}_0 \overline{m} w(r) = 0 & \text{for } 0<r<1 \; , \\
w''(r) + \frac{N-1}{r}w'(r)-\tilde{\lambda}_0 \underline{m} w(r) = 0 & \text{for } r>1 \; , \\
w'(0) = 0 \; , \\
\lim_{r \to +\infty} w(r) = 0 \; , \\
w(r) \in C^1(\R^{+}) \; .
\end{cases}
\end{equation}
In particular, 
\begin{equation}
w(r) =
\begin{cases}\label{eq:SolRnShape}
A \, r^{1-N/2} J_{N/2-1} \left(\sqrt{\tilde{\lambda}_0 \overline{m}}r \right) & \text{for } 0\le r <1 \; , \\
B \, r^{1-N/2} K_{N/2-1} \left(\sqrt{\tilde{\lambda}_0 \underline{m}}r \right) & \text{for } r \ge 1 \; , \\
\end{cases}
\end{equation}
where $J_\nu$, $K_\nu$ denote the the Bessel function of the first kind and the special one of the second kind, and $A$ and $B$ are positive suitable constants.

Moreover,
\[
w(r) \sim \frac{C}{\sqrt{r}}\, e^{-r\sqrt{\tilde{\lambda}_0 \underline{m}}} \quad \text{for } r \to +\infty, \, C>0 \; ,
\]
and
\[
\frac{w(r)}{w'(r)} \sim -\sqrt{\tilde{\lambda}_0 \underline{m}} \quad \text{for } r \to +\infty \; .
\]
\end{proposition}

System \eqref{eq:ODEsystem1} is simply equation \eqref{eqn:sulutionRn} written with respect to the radial coordinate, exploiting the regularity of $w$. Solution \eqref{eq:SolRnShape} comes from the well established properties of the Bessel functions, see for instance 
\cite[Ch. 9]{AbramowitzStegun:Handbook}.

\section{Blow-up procedure}\label{blowUp}

This section marks the beginning of the study of the asympotics for the first generalized eigenvalue with spherical weights $\lambda_{\varepsilon}$, i.e the study of problem \eqref{sphericalProblem}. As we mentioned, we do this through a blow-up procedure.

According to Proposition \ref{prop:PropositionMinimum} we denote with 
$u_{\varepsilon}$ the family of positive eigenfunctions associated to 
$\lambda_{\varepsilon}$, normalized in $L^2(\Omega)$, and $\mathbf{x}_{\varepsilon} \in \hat{\Omega}_{\varepsilon}$ the center of the sphere for which $\lambda_{\varepsilon}$ is attained. Since $\Omega$ is bounded, from any sequence $\varepsilon_n \to 0$ we can extract a subsequence, still denoted by $\varepsilon_n$, such that $\mathbf{x}_{\varepsilon_n} \to \bar{\mathbf{x}} \in \overline{\Omega}$ for $n \to +\infty$. 
Finally, we recall that the blow-up sequences $\tilde u_n= \tilde u_{\eps_n}$, 
$\tilde \Omega_n= \tilde \Omega_{\eps_n}$, $\tilde \lambda_n= \tilde \lambda_{\eps_n}$
defined in equations \eqref{eq:kbeta_intro}, \eqref{eq:blowupseq_intro}, \eqref{eq:lambdabu_intro}, satisfy problem \eqref{GeneralDiffProblemBU_intro}. We remark that, even in the blow up scale, $\tilde{\lambda}_{n}$ satisfies:
\begin{equation}\label{eq:bu_min}
\tilde{\lambda}_{n} = \inf_{B^1(\mathbf{x}) \subset \tilde{\Omega}_n} \lambda^1(B^1(\mathbf{x}), \tilde{\Omega}_n) \; ,
\end{equation}
which is actually a minimum and is attained by $\tilde{u}_n$ and $\mathbf{x} = \mathbf{0}$. This yields the convergence of the sequence $\tilde{\lambda}_{n}$.
\begin{lemma}\label{lem:conv2la_0}
Let $\tilde \lambda_0>0$ be defined as in Proposition \ref{result:ResultProblemOverRn}. Then
\[
\lim_{n \to +\infty} \tilde{\lambda}_{n} = \tilde{\lambda}_{0}.
\] 
\end{lemma}
\begin{proof}
On the one hand, we can use (the trivial extension of) $\tilde{u}_n\in H^1(\R^N)$ and 
$\tilde m_0$ in \eqref{lambda0defRn} to infer
\[
\tilde{\lambda}_{0} \le \tilde{\lambda}_{n},\qquad\text{for every }n.
\]
On the other hand, take any $\hat y \in \Omega$, with $d(\hat \by,\partial \Omega)=
\delta>0$ and define $\by_n=(\hat\by-\bx_n)/k_n$. Then if $n$ is sufficiently large 
we obtain  $B^1(\mathbf{y}_n)\subset B^{\delta/k_n}(\mathbf{y}_n)  \subset 
\tilde{\Omega}_n$, and \eqref{eq:bu_min} yields
\[
\tilde{\lambda}_{n} \le  \lambda^1(B^1(\mathbf{y}_n),B^{\delta/k_n}(\mathbf{y}_n)) =\lambda^1(B^1(\mathbf{0}),B^{\delta/k_n}(\mathbf{0})).
\] 
Since $\delta/k_n$ diverges we can conclude by using 
\eqref{eqn:characterizationTildeLambda0}.
\end{proof}

Now we present a result which guarantees that the $\mathbf{x}_{\varepsilon_n}$ cannot reach $\partial\Omega$ too fast:
\begin{proposition}\label{prop:PropAwayFromBoundary}
\[
d_n \coloneqq d(\mathbf{x}_{\varepsilon_n}, \partial\Omega)/k_n \to +\infty \qquad \text{for } n \to +\infty \; .
\]
In particular, for any $K \subset \subset \R^N$, $K \subset \tilde{\Omega}_n$ for $n$ big enough.
\end{proposition}
We postpone the proof of it after Proposition \ref{PropC1alphaConvergence}.

Assuming for the moment Proposition \ref{prop:PropAwayFromBoundary}, we are ready to prove the first convergence result for $\tilde{u}_n$:
\begin{proposition}\label{PropC1alphaConvergence}
There exists $\bar{u} \in H^1(\R^N) $ such that: 
\begin{itemize}
\item[(i)] $ \tilde{u}_n \to \bar{u} $ weakly in $ H^1(\R^N) $ for $ n \to +\infty $, up to a subsequence that we still denote $\tilde{u}_n$,
\item[(ii)] for any fixed $0 < \alpha < 1$, there exists a further subsequence $\tilde{u}_{n_k}$ such that
\[
\tilde{u}_{n_k} \to \bar{u} \;\text{in} \,\, C^{1, \alpha}(K) \quad \text{if } n \to +\infty, \; \forall \, K \subset \subset \R^N \; .
\]
In particular, $\bar{u} \in C^{1, \alpha}(K)$ for any $0 < \alpha < 1$, $K \subset \subset \R^N$ and $\bar{u} \in C^{1}(\R^N)$,
\item[(iii)] $\bar{u}$ satisfies:
\[
\int_{\R^n} \nabla \bar{u} \cdot \nabla v = \tilde{\lambda}_0 \int_{\R^n} \tilde{m}_0 \bar{u} v \qquad \forall v \in H^1(\R^N) \; ,
\]
where $\tilde{\lambda}_0$, $\tilde{m}_0$ are as in Proposition \ref{prop:PropositionBlowUp},
\item[(iv)] $\bar{u} > 0$ in $\R^N$.
\end{itemize}
\end{proposition}
\begin{proof}
(i) We consider $\tilde{u}_n$ extended to 0 outside of $\tilde{\Omega}_n$. Now, since $\| \tilde{u}_n \|_{L^2(\R^N)} = 1$ and since they solve \eqref{GeneralDiffProblemBU_intro}, considering as test function $\tilde{u}_n$ in its weak formulation, we get:
\begin{equation}\label{eq:PropC1alphaConvergence1}
\| \nabla \tilde{u}_n \|^2_{L^2(\R^N)} \le \tilde{\lambda}_{n} \| \tilde{m}_0 \|_{L^{\infty}(\R^N)} \| \tilde{u}_n \|_{L^2(\R^N)} \le C \; ,
\end{equation}
thanks to Lemma \ref{lem:conv2la_0} and the definition of $\tilde{m}_0$. The constant $C$ does not depend on $n$. This implies that $\tilde{u}_n$ converges weakly in $H^1(\R^N)$ to a function $\bar{u} \in H^1(\R^N)$, up to a subsequence that we still denote $\tilde{u}_n$. 

(ii) By standard elliptic regularity theory, since $\tilde{m}_0 \in L^{\infty}(\R^N)$, we have that $\tilde{u}_n \in C^{1, \alpha}(\tilde{\Omega}_n)$ for any $0 < \alpha < 1$. Actually, since the $H^1(\tilde{\Omega}_n)$ norm of $\tilde{u}_n$ is bounded uniformly with respect to $n$, also their $C^{1, \alpha}(K)$ norm is bounded uniformly in $n$, for any compact set $K \subset \R^N$ and for $n$ big enough, since by Proposition \ref{prop:PropAwayFromBoundary} we have $K \subset \tilde{\Omega}_n$. Thus, since $C^{1, \beta}(K) \hookrightarrow \hookrightarrow C^{1, \alpha}(K) $ for any $0 < \alpha <  \beta < 1$, $K \subset \subset \R^N$, we can fix $\alpha$ and extract a subsequence $\tilde{u}_{n_k}$ which converges in $C^{1, \alpha}(K) $ to some $f \in C^{1, \alpha}(K)$, which must coincide with $\bar{u}$ due to point (i).

Actually, choosing $K$ as a sequence of closed balls centered at the origin with divergent radius and proceeding with a diagonal argument,  for fixed $0 < \alpha < 1$, we can find a subsequence $\tilde{u}_{n_k}$ such that:
\[
\tilde{u}_{n_k} \to \bar{u} \;\text{in} \,\, C^{1, \alpha}(K) \quad \text{if } n \to +\infty, \; \forall \, K \subset \subset R^N \; .
\]
This, in particular, implies that:
\[
\bar{u} \in C^{1, \alpha}(K) \quad \forall \, 0 < \alpha < 1, \, K \subset \subset \R^N \; ,
\]
and also
\[
\bar{u} \in C^{1}(\R^N) \; .
\]
We also remark that the previous argument tells us that
\begin{equation}\label{eqn:nonNegativeness}
\bar{u} \ge 0 \quad \text{in } \R^N \; ,
\end{equation}
since over any compact set it is the uniform limit of the sequence of positive functions $\tilde{u}_{n_k}$. 

(iii)  Choose an arbitrary $v \in C_c^{\infty}(\R^N)$. Then, thanks to Proposition \ref{prop:PropAwayFromBoundary} and for $n$ big enough, $\supp(v) \subset \tilde{\Omega}_n$. Thus $\tilde{u}_n$ satisfies:
\[
\int_{\R^N} \nabla \tilde{u}_n \cdot \nabla v = \tilde{\lambda}_n \int_{\R^N} \tilde{m}_0 \tilde{u}_n v  \; .
\]
Now, using Lemma \ref{lem:conv2la_0} and point (i), we can pass to the limit 
obtaining
\[
\int_{\R^N} \nabla \bar{u} \cdot \nabla v = \tilde{\lambda}_0 \int_{\R^N} \tilde{m}_0 \bar{u} v  \; .
\]
By the density of $C_c^{\infty}(\R^N)$ in $H^1(\R^N)$ and by the arbitrariness of $v$, the above relation actually holds for any $v \in H^1(\R^N)$.

(iv) We recall that by hypothesis, $\int_{\R^N} \tilde{m}_n \tilde{u}_n^2 > 0$ for any $n \in \N$. This relation can be rewritten as: 
\begin{equation}\label{eqnd:positivenessEstimate1}
\overline{m}\int_{B^1(\mathbf{0})} \tilde{u}_n^2 - \underline{m}\int_{B^1(\mathbf{0})^c} \tilde{u}_n^2 > 0 \; .
\end{equation}
On the other hand, we also have that $\| \tilde{u}_n \|_{L^2(\R^N)} = 1$, that we rewrite as:
\begin{equation}\label{eqnd:positivenessEstimate2}
\underline{m} \int_{B^1(\mathbf{0})} \tilde{u}_n^2 + \underline{m} \int_{B^1(\mathbf{0})^c} \tilde{u}_n^2 = \underline{m} \; .
\end{equation}
Summing up \eqref{eqnd:positivenessEstimate1} and \eqref{eqnd:positivenessEstimate2} we obtain:
\begin{equation}\label{eqnd:positivenessEstimate3}
\int_{B_1(\mathbf{0})} \tilde{u}_n^2 > \frac{\underline{m}}{\underline{m} + \overline{m}} \; .
\end{equation}
On the other hand, since $\mathcal{L}(B^1(\mathbf{0})) = 1$, we also have:
\begin{equation}\label{eqnd:positivenessEstimate10}
\int_{B^1(\mathbf{0})} \tilde{u}_n^2 \le \max_{\overline{B}^1(\mathbf{0})} \tilde{u}_n^2 \; .
\end{equation}
Thus, combining \eqref{eqnd:positivenessEstimate3} and \eqref{eqnd:positivenessEstimate10} we have:
\begin{equation}\label{eqnd:positivenessEstimate4}
\max_{\overline{B}^1(\mathbf{0})} \tilde{u}_n \ge \left( \frac{\underline{m}}{\underline{m} + \overline{m}} \right)^{\frac{1}{2}} \; ,
\end{equation}
where the right hand side is strictly positive and independent on $n$. Thus, from the uniform convergence obtained in point (ii), we have:
\begin{equation}\label{eqnd:positivenessEstimate5}
\max_{\overline{B}^1(\mathbf{0})} \bar{u} > 0 \; .
\end{equation}
Now, since by point (iii) $\bar{u}$ satisfies
\[
-\Delta \bar{u} + \tilde{\lambda}_0 \underline{m} \bar{u} = \tilde{\lambda}_0 (\underline{m} + \overline{m}) \ind{B_1(\mathbf{0})}  \bar{u} \quad \text{in } \R^N \; , 
\]
where the right hand side is non-negative thanks to \eqref{eqn:nonNegativeness}, using the Strong Maximum Principle we have that either $\bar{u} > 0$, or $\bar{u} \equiv 0$ in $\R^N$. To conclude, thanks to \eqref{eqnd:positivenessEstimate5} we must have $\bar{u} > 0$ in $\R^N$.
\end{proof}

Now we turn to the proof of Proposition \ref{prop:PropAwayFromBoundary}.
\begin{proof}[Proof of Proposition \ref{prop:PropAwayFromBoundary}]
We proceed by contradiction. Suppose that $d_n$ is bounded. Then, up to a subsequence that we denote with the same index $n$, we have $\mathbf{x}_{\varepsilon_n} \to \mathbf{x}' \in \partial\Omega$ for $n \to +\infty$. This implies that, thanks to the smoothness of $\Omega$ and proceeding similarly to Proposition \ref{PropC1alphaConvergence} (i), (ii), we can find an open half-space $T \subset \R^N$ and $\bar{u} \in H_0^1(T)$ which is the weak limit of $\tilde{u}_n$ in $H^1(\R^N)$ (extended to $0$ outside of $\tilde{\Omega}_n$), such that:
\begin{equation}\label{eqn:ContradictionDistanceBoundaryEqn}
\int_{\R^n} \nabla \bar{u} \cdot \nabla v = \tilde{\lambda}_0 \int_{\R^n} \tilde{m}_0 \bar{u} v \qquad \forall v \in H_0^1(T) \; .
\end{equation}
Moreover, proceeding similarly as in the proof of Proposition \ref{PropC1alphaConvergence} (iii), (iv), we have that for any fixed $0 < \alpha < 1$, there exists a further subsequence $\tilde{u}_{n_k}$ such that $\tilde{u}_{n_k} \to \bar{u}$ in $C^{1, \alpha}(K)$ for $n \to +\infty$ and for any $K \subset \subset T$, and $\bar{u} > 0$ in $T$ by the Strong Maximum Principle. 

We also remark that we need to have $\int_{T} \tilde{m}_0 \bar{u}^2 > 0$. Otherwise, choosing $v = \bar{u}$ in \eqref{eqn:ContradictionDistanceBoundaryEqn} we woud have $\int_{T}|\nabla \bar{u}|^2 = 0$ and thus $\bar{u}$ constant on $T$. Since $\bar{u} \in H_0^1(T)$ this constant would be zero and hence $\bar{u} = 0$ on $T$, against the fact that $\bar{u} > 0$ on $T$. 

Thus, from the above discussion we can write:
\begin{equation}\label{eqn:rayleighQuotientContradictionDistanceBoundary}
\tilde{\lambda}_0 = \frac{\int_{T}|\nabla \bar{u}|^2}{\int_{T} \tilde{m}_0 \bar{u}^2} = \frac{\int_{\R^N}|\nabla \bar{u}|^2}{\int_{\R^N} \tilde{m}_0 \bar{u}^2} \; .
\end{equation}
Now, thanks to Proposition \ref{result:ResultProblemOverRn} we have that $\bar{u}$ needs to be radially symmetric with respect to a point $\mathbf{x} \in \R^N$. However, since $\bar{u} = 0$ on $\partial T$ we necessarily have that $\bar{u} = 0$. But this is an absurd, since $\bar{u} > 0$ in $T$.

Thus, we have reached a contradiction and so $d_n \to +\infty$ for $n \to + \infty$.
\end{proof}
Now we proceed to a first identification of the function $\bar{u}$.
\begin{proposition}\label{prop:PropLimitToProblemOverRn}
The following results hold: 
\[
\bar{u} = A w\qquad\text{ for some positive constant $A$}.
\]
\end{proposition}
\begin{proof}
To begin with, we remark that $\int_{\R^N} \tilde{m}_0 \bar{u}^2 > 0$. Indeed, otherwise by Proposition \ref{PropC1alphaConvergence} (iii) we would have $\int_{\R^N}|\nabla \bar{u}|^2 = 0$ and thus $\bar{u} = c$, with $c \in \R$. But, since $\bar{u} \in H^1(\R^N)$ we would have $\bar{u} = 0$. However, this contradicts Proposition \ref{PropC1alphaConvergence} (iv). Thus, we can write:
\[
\tilde{\lambda}_0 = \frac{\int_{\R^N}|\nabla \bar{u}|^2}{\int_{\R^N} \tilde{m}_0 \bar{u}^2} \; .
\]
We conclude by the uniqueness stated in Proposition \ref{result:ResultProblemOverRn} 
and by the positiveness of $\bar{u}$, stated in Proposition 
\ref{PropC1alphaConvergence} (iv).
\end{proof}

\section{Convergence in \texorpdfstring{$H^1(\R^N)$}{H1(RN)}}\label{section:H1Convergence} 
This section is devoted to the conclusion of the proof of Proposition \ref{prop:PropositionBlowUp}. What remains to be proved is that $A = 1$ in Proposition \ref{prop:PropLimitToProblemOverRn} and that, up to a further subsequence, the convergence to $w$ in $H^1(\R^N)$ is strong. We recall that, unless otherwise stated, we choose $\tilde{u}_n$ and $w$ normalized in $L^2(\R^N)$.

Our arguments will be based on the well known concentration-compactness argument by P. L. Lions \cite[Lemma I.1]{Lions:concentrCompact}, that we reformulate in the following way
\begin{lemma}\label{result:LemmaCompactVanishDichotomy}
Let $\rho_n$ be a sequence in $L^1(\R^N)$ satisfying:
\[
\rho_n \ge 0 \; \text{in } \R^N, \quad \int_{\R^N} \rho_n = \lambda
\]
where $\lambda > 0$ is fixed. Then, there exists a subsequence $\rho_{n_k}$ satisfying one of the three following possibilities:
\begin{itemize}
\item[(i)] (compactness) there exists $\mathbf{y}_k \in \R^N$ such that:
\[
\forall \, \varepsilon > 0, \; \exists \, R < +\infty \; : \; \int_{\mathbf{y}_k + B_R(\mathbf{y}_k)} \rho_{n_k} \ge \lambda - \varepsilon \; ,
\]
\item[(ii)] (vanishing) 
\begin{equation}\label{eqn:vanishingDef}
\forall \, R>0 \; : \; \lim_{k \to +\infty} \sup_{\mathbf{y} \in \R^N} \int_{\mathbf{y} + B_R(\mathbf{y})} \rho_{n_k} = 0  \; ,
\end{equation}
\item[(iii)] (dichotomy) there exists $\alpha \in (0, \lambda)$ such that, for all $\varepsilon > 0$ there exist $\mathbf{y}_k \in \R^N$, $R > 0$, $R_k$ with $R_k \to +\infty$ for $k \to +\infty$ and $k_0 \ge 1$ such that for $k \ge k_0$:
\begin{equation}\label{eqn:dichotomyDef1}
\int_{R \le |\mathbf{x} - \mathbf{y}_k | \le R_k} \rho_{n_k} \le \varepsilon \; ,
\end{equation}
\begin{equation}\label{eqn:dichotomyDef2}
\left\vert \int_{\mathbf{y}_k + B_R(\mathbf{y}_k)} \rho_{n_k} - \alpha \right\vert \le \varepsilon \; ,
\end{equation}
\begin{equation}\label{eqn:dichotomyDef3}
\left\vert \int_{\R^N \setminus (\mathbf{y}_k + B_{R_k}(\mathbf{y}_k))} \rho_{n_k} - (\lambda - \alpha) \right\vert \le \varepsilon \; ,
\end{equation}
\end{itemize}
\end{lemma}
To fix notations, in our case $\rho_n = \tilde{u}_n^2$ and $\lambda = 1$. 

We will prove that for the blow-up sequence $\tilde{u}_n$ introduced in Section \ref{blowUp}, only the compactness case can occur. More precisely, our aim is to prove the following
\begin{proposition}\label{prop:PropCompactnessProperty}
Up to a further subsequence that we still denote $\tilde{u}_n$,
\[
\forall \, \varepsilon > 0, \; \exists \, R < +\infty \; : \; \int_{B_R(\mathbf{0})} \tilde{u}_n^2 \ge 1 - \varepsilon \; .
\]
\end{proposition}
We proceed by excluding both vanishing and dichotomy. This is done in the following lemmas.
\begin{lemma}\label{lemma:LemmaNoVanishing}
For any subsequence, that we still denote $\tilde{u}_n$, vanishing is not possible.
\end{lemma}
\begin{proof}
Consider $\mathbf{y} = 0$ and choose $R$ such that $B_R(\mathbf{0}) = B^1(\mathbf{0})$ . Then, 
\[
\sup_{\mathbf{y} \in \R^N} \int_{\mathbf{y} + B_R(\mathbf{y})} \tilde{u}_{n} \ge \int_{B_R(\mathbf{0})} \tilde{u}_{n} \ge \frac{\underline{m}}{\underline{m} + \overline{m}}  \; ,
\]
uniformly in $n$. The last inequality comes from \eqref{eqnd:positivenessEstimate3}. The above relation implies that vanishing cannot happen.
\end{proof}

Now we can state and prove the following
\begin{lemma}\label{lemma:LemmaNoDichotomy}
For any subsequence, that we still denote $\tilde{u}_n$, dichotomy is not possible.
\end{lemma}
\begin{proof}
In the following, $C$ will denote a fixed constant. 

We begin by noticing that, for any $R>0$ big enough and choosing $n$ big enough, on $\R^N \setminus B_R(\mathbf{0})$ we have:
\begin{equation}\label{eqn:estimateUn}
\tilde{u}_n \le C e^{-\sqrt{\tilde{\lambda}_0 \underline{m}}R} \; ,
\end{equation}
where the constant $C$ is independent of $R$. This is a direct consequence of the Maximum Principle, of Proposition \ref{PropC1alphaConvergence} (iii) and of Proposition \ref{result:ResultShapeSolutionOverRn}.

We also remark that, for any fixed $R > 0$, in Lemma \ref{result:LemmaCompactVanishDichotomy} (iii) we can consider $\mathbf{y}_n = \mathbf{0}$. Indeed, it is sufficient to prove that $\mathbf{y}_n$ is bounded. Suppose it is not. Consider a divergent subsequence $\mathbf{y}_{n_k}$. Then, applying a diagonal argument, thanks to \eqref{eqn:estimateUn} for any $\varepsilon' > 0$ we can extract a subsequence $\tilde{u}_{n_k}$ such that $\int_{\mathbf{y}_{n_k} + B_R(\mathbf{y}_{n_k})} \tilde{u}_{n_k}^2 \le \varepsilon'$. But this contradicts \eqref{eqn:dichotomyDef2}. Thus, from now on, we will consider $\mathbf{y}_n = \mathbf{0}$.

Now, consider a sequence $\varepsilon_k \to 0^{+}$. Consequently, in Lemma \ref{result:LemmaCompactVanishDichotomy} (iii), in place of the radius $R$ we can  consider a sequence $R_k \to +\infty$. We define a family of cutoff functions for the sets $B_{R_k}(\mathbf{0})$, as
\[
\eta_k \coloneqq
\begin{cases}
1 & \text{in } \overline{B}_{R_k}(\mathbf{0}) \; ,\\
1 - \bigl( |\mathbf{x}| - R_k)^2 & \text{in } B_{R_k+1}(\mathbf{0}) \setminus B_{R_k}(\mathbf{0}) ,\\
0 & \text{otherwise} .
\end{cases}
\]
We remark that $0 \le \eta_k \le 1$ and that the functions $\tilde{u}_n\eta_k$ belong to $H_0^1(\R^N)$.

Now, for any fixed $\varepsilon_k$ small enough, for $n$ big enough
\begin{equation}\label{eqn:dichotomy0}
\int_{R_k \le |\mathbf{x}| \le R_k+1} \tilde{u}_n^2 \eta_k^2 \le C e^{-\sqrt{\tilde{\lambda}_0 \underline{m}}R_k} R_k = o(1) \quad \text{for } k \to +\infty  \; ,
\end{equation}
thanks to \eqref{eqn:estimateUn}. Moreover, by direct computation
\[
\int_{R_k \le |\mathbf{x}| \le R_k+1} | \nabla \eta_k |^2 \le C R_k
\]
and thus for $n$ big enough
\begin{equation}\label{eqn:dichotomy1}
\int_{R_k \le |\mathbf{x}| \le R_k+1} \tilde{u}_n^2 | \nabla \eta_k |^2 \le C e^{-\sqrt{\tilde{\lambda}_0 \underline{m}}R_k} R_k  \; ,
\end{equation}
\begin{equation}\label{eqn:dichotomy2}
\left\vert 2 \int_{R_k \le |\mathbf{x}| \le R_k+1} \tilde{u}_n \eta_k \nabla \tilde{u}_n \cdot  \nabla \eta_k \right\vert \le C e^{-\sqrt{\tilde{\lambda}_0 \underline{m}}R_k} R_k  \; ,
\end{equation}
and thanks to Proposition \ref{result:ResultShapeSolutionOverRn} and Proposition \ref{PropC1alphaConvergence} (iii), for any $\varepsilon' > 0$ there exists $k$ such that, for $n$ big enough:
\begin{equation}\label{eqn:dichotomy3}
\int_{R_k \le |\mathbf{x}| \le R_k+1} \eta_k^2 | \nabla \tilde{u}_n |^2 \le \varepsilon'  \; .
\end{equation}
Hence, from \eqref{eqn:dichotomy1}, \eqref{eqn:dichotomy2} and \eqref{eqn:dichotomy3}, up to the extraction of a diagonal subsequence $\tilde{u}_{n_k}$ (where now $n$ depends on $k$)
\begin{equation}\label{eqn:dichotomyFirst}
\int_{R_k \le |\mathbf{x}| \le R_k+1} | \nabla (\tilde{u}_{n_k} \eta_k) |^2 \le o(1) \quad \text{for } k \to +\infty  \; ,
\end{equation}
and analogously, from \eqref{eqn:dichotomy0}
\begin{equation}\label{eqn:dichotomySecond}
\int_{R_k \le |\mathbf{x}| \le R_k+1} \tilde{u}_{n_k}^2 \eta_k^2 \le o(1) \quad \text{for } k \to +\infty  \; .
\end{equation}
Now we consider $\tilde{u}_{n_k}$ normalized with $\| \nabla \tilde{u}_{n_k} \|_{L^2(\R^N)} = 1$. If we denote with $c_{n_k}$ the normalization constants, thanks to Proposition \ref{PropC1alphaConvergence} we have $0 < C \le c_{n_k} \le M < +\infty$, with $C$ and $M$ positive constants. Moreover, without loss of generality, in our subsequence $\tilde{u}_{n_k}$ we may assume that $n_k > k_0$  for any $k \ge 0$, where $k_0$ is defined in Lemma \ref{result:LemmaCompactVanishDichotomy} (iii). In this way we ensure that \eqref{eqn:dichotomyDef1}, \eqref{eqn:dichotomyDef2} and \eqref{eqn:dichotomyDef3} hold, with $c_{n_k}^2 \alpha$ and $c_{n_k}^2 (\lambda - \alpha)$ in place of $\alpha$ and $\lambda - \alpha$ respectively. 

To conclude, thanks to \eqref{eqn:dichotomyFirst}, \eqref{eqn:dichotomySecond}, \eqref{eqn:dichotomyDef1}, \eqref{eqn:dichotomyDef2} and \eqref{eqn:dichotomyDef3}, we have:
\begin{align*}
& \frac{1}{\tilde{\lambda}_{n_k}} = \int_{\R^N} \tilde{m}_{n_k} \tilde{u}_{n_k}^2 = \int_{B_{R_k + 1}(\mathbf{0})} \tilde{m}_{n_k} \tilde{u}_{n_k}^2 \eta_k^2 + o(1) + \int_{R_k \le |\mathbf{x}| \le R_{n_k}} \tilde{m}_{n_k} \tilde{u}_{n_k}^2  +\\
& + \int_{B_{R_{n_k}}^c(\mathbf{0})} \tilde{m}_{n_k} \tilde{u}_{n_k}^2 \le  \frac{1}{\tilde{\lambda}_{n_k}} \int_{B_{R_k+1}(\mathbf{0})} | \nabla (\tilde{u}_{n_k} \eta_k) |^2 - \underline{m} \int_{B_{R_{n_k}}^c(\mathbf{0})} \tilde{u}_{n_k}^2 + \\
& + o(1) + \varepsilon_k = \frac{1}{\tilde{\lambda}_{n_k}} \int_{B_{R_k}(\mathbf{0})} | \nabla \tilde{u}_{n_k} |^2 + \frac{1}{\tilde{\lambda}_{n_k}} \int_{R_k \le |\mathbf{x}| \le R_k+1} | \nabla (\tilde{u}_{n_k} \eta_k) |^2 + 2 \varepsilon_k + \\
& + o(1) - \underline{m} c_{n_k}^2 (\lambda-\alpha) \le \frac{1}{\tilde{\lambda}_{n_k}}  - \underline{m} c_{n_k}^2 (\lambda-\alpha) + 2o(1) + 2 \varepsilon_k < \frac{1}{\tilde{\lambda}_{n_k}} 
\end{align*}
for $k$ big enough, which is an absurd. This concludes the proof.
\end{proof}
\begin{proof}[Proof of Proposition \ref{prop:PropCompactnessProperty}]
Combining Lemma \ref{result:LemmaCompactVanishDichotomy} and Lemmas \ref{lemma:LemmaNoVanishing}, \ref{lemma:LemmaNoDichotomy}, the only thing left to prove is that we can choose $\mathbf{y}_n = \mathbf{0}$ for any $n \in \N$. Actually, it is sufficient to prove that $\mathbf{y}_n$ is bounded in $\R^N$.

Suppose it is not. Then we choose $\varepsilon$ and extract an unbounded subsequence $\mathbf{y}_{n_k}$ such that, for $R$ big enough independent of $k$ and $k$ big enough:
\[
\mathbf{y}_{n_k} \cap B_1(\mathbf{0}) = \emptyset, \quad \varepsilon < \frac{\underline{m}}{\underline{m} + \overline{m}} \; .
\]
Thanks to \eqref{eqnd:positivenessEstimate3}, we have reached an absurd and hence $\mathbf{y}_n$ is bounded.
\end{proof}
We are ready to prove the $L^2(\R^N)$ convergence of $\tilde{u}_n$ to $w$:
\begin{proposition}\label{prop:L2ConvergenceRN}
If $\| \tilde{u}_n \|_{L^2(\R^N)} = 1$ and $\| w \|_{L^2(\R^N)} = 1$, then
\[
\tilde{u}_n \to w \quad \text{in } L^2(\R^N) \; \text{for } n \to +\infty 
\]
\end{proposition}
\begin{proof}
Choose $\varepsilon > 0$. Using Proposition \ref{prop:PropCompactnessProperty} we can choose $R > 0$ such that $ \int_{B_R^c(\mathbf{0})} \tilde{u}_n^2 \le \varepsilon  $. Moreover, since $w \in L^2(\R^N)$, we can choose $R > 0$ such that $ \int_{B_R^c(\mathbf{0})} w^2 \le \varepsilon$. Thus:
\begin{align*}
\int_{\R^N} |\tilde{u}_n - w |^2 = & \int_{B_R(\mathbf{0})} |\tilde{u}_n - w |^2 + \int_{B_R^c(\mathbf{0})} |\tilde{u}_n - w |^2 \le \\
& \int_{B_R(\mathbf{0})} |\tilde{u}_n - w |^2 + 2\varepsilon \le 3\varepsilon \; ,
\end{align*}
for $n$ big enough, where in the last inequality we have used Proposition \ref{PropC1alphaConvergence} (iii). This concludes the proof.
\end{proof}
We remark that, exploiting exponential estimates in a way similar to \cite[Prop. 3.4]{NiWei:Spike95}, one can actually prove convergence of $\tilde{u}_n$ to $w$ in $L^p(\R^N)$ for any $1 \le p \le +\infty$. 
\begin{proof}[End of the proof of Proposition \ref{prop:PropositionBlowUp}]
We are only left to prove that $\| \nabla \tilde{u}_n - \nabla w \|_{L^2(\R^N)} \to 0$ for $n \to +\infty$. Actually, since $\tilde{u}_n \rightharpoonup w$ in $H^1(\R^N)$, we only need to prove that
\begin{equation}\label{eqn:GradNormConvergence}
\| \nabla \tilde{u}_n \|_{L^2(\R^N)}^2 \to \| \nabla w \|_{L^2(\R^N)}^2 \quad \text{for } n \to +\infty \; .
\end{equation}
Testing the weak forms of \eqref{problemOverRn} and \eqref{GeneralDiffProblemBU_intro} with $w$ and $\tilde{u}_n$ respectively, we get:
\begin{align*}
& | \| \nabla \tilde{u}_n \|_{L^2(\R^N)}^2 - \| \nabla w \|_{L^2(\R^N)}^2  | \le \\
& \le \left\vert (\tilde{\lambda}_n - \tilde{\lambda}_0) \int_{\R^N} \tilde{m}_0 \tilde{u}_n^2 \right\vert + \left\vert \tilde{\lambda}_0 \int_{\R^N} \tilde{m}_0 (\tilde{u}_n^2 - w^2) \right\vert \le \\
& \le (\tilde{\lambda}_n - \tilde{\lambda}_0) \| \tilde{m}_0 \|_{L^{\infty}(\R^N)} \| \tilde{u}_n \|_{L^2(\R^N)}^2 + \\
& + \| \tilde{m}_0 \|_{L^{\infty}(\R^N)} \| \tilde{u}_n + w \|_{L^2(\R^N)}^2 \| \tilde{u}_n - w \|_{L^2(\R^N)}^2 \to 0
\end{align*}
for $n \to +\infty$, since $\tilde{u}_n$ and $w$ are bounded in $L^2(\R^N)$, $\tilde{\lambda}_n \to \tilde{\lambda}_0$ and thanks to Proposition \ref{prop:L2ConvergenceRN}. This concludes the proof.
\end{proof}

\section{Estimates from above}\label{section:EstimateAbove}
Thanks to Proposition \ref{prop:PropositionBlowUp} up to a subsequence we have qualitative results of convergence of $\tilde{u}_{\varepsilon}$ and $\tilde{\lambda}_{\varepsilon}$. To proceed, we wish to look for quantitative results, i.e. velocity of convergence of both the sequence of blow-up solutions and of the eigenvalues, and the positioning of the blow-up points $\mathbf{x}_{\varepsilon}$ in ${\Omega}$. 

The answers to these questions are contained in Theorem \ref{theorem:TheoremAsymptExpansion}, that we proceed to prove in this section and in the next one. In this section, in particular, we present two estimates from above for the eigenvalue $\lambda_{\varepsilon}$.

We denote with $\mathbf{q} \in \Omega$ a point at maximum distance from 
$\partial\Omega$, and
\[
\begin{split}
\tilde{\Omega}_{\varepsilon} &\coloneqq \{ \mathbf{x} \in \R^N : k_{\varepsilon} \mathbf{x} + \mathbf{x}_{\varepsilon} \in \Omega \} \; ,\\
\Omega_{\varepsilon} &\coloneqq \{ \mathbf{x} \in \R^N : k_{\varepsilon} \mathbf{x} + \mathbf{q} \in \Omega \} \; .
\end{split}
\]
where, as usual, $k_{\varepsilon} \coloneqq \varepsilon^{{1}/{N}}$. 
Then, following \cite{NiWei:Spike95}, we consider the $H^1_0(\Omega_{\varepsilon})$ projection of $w$, denoted $P_{\Omega_{\varepsilon}} w$ and the $H^1_0(\tilde{\Omega}_{\varepsilon})$ projection of $w$, denoted $P_{\tilde{\Omega}_{n}} w$. 

In particular, defining 
\[
f_0(w) \coloneqq \tilde{\lambda}_0 (\overline{m}+\underline{m}) \ind{B^1(\mathbf{0})} w
\] 
the function $P_{\Omega_{\varepsilon}} w$ solves
\begin{equation}
\begin{cases}\label{eq:Pw}
-\Delta P_{\Omega_{\varepsilon}} w  + \tilde{\lambda}_0 \underline{m} P_{\Omega_{\varepsilon}} w = f_0(w) & \text{in } \Omega_{\varepsilon} \; , \\
P_{\Omega_{\varepsilon}} w = 0 & \text{on } \partial\Omega_{\varepsilon} \; ,
\end{cases}
\end{equation}
while $P_{\tilde{\Omega}_{\varepsilon}} w$ solves
\begin{equation}
\begin{cases}\label{eq:Ptildew}
-\Delta P_{\tilde{\Omega}_{\varepsilon}} w + \tilde{\lambda}_0 \underline{m} P_{\tilde{\Omega}_{\varepsilon}} w = f_0(w) & \text{in } \tilde{\Omega}_{\varepsilon} \; , \\
P_{\tilde{\Omega}_{\varepsilon}} w = 0 & \text{on } \partial\tilde{\Omega}_{\varepsilon} \; .
\end{cases}
\end{equation}
Notice that, by the Strong Maximum Principle, both $P_{\Omega_{\varepsilon}} w$ and $P_{\tilde{\Omega}_{\varepsilon}} w$ are positive on their domain of definition.

If we introduce the quantity $\beta_{\varepsilon} \coloneqq 1/k_{\varepsilon}$, and the functions
\begin{equation}\label{PsiEpsilonDefinition}
\Psi_{\varepsilon}(\mathbf{x}) \coloneqq -k_{\varepsilon} \log(w - P_{\Omega_{\varepsilon}} w)\left(\frac{\mathbf{x}-\mathbf{q}}{k_{\varepsilon}}\right) \quad \text{on } \Omega \; ,
\end{equation}
and
\begin{equation}\label{PsiTildeEpsilonDefinition}
\tilde{\Psi}_{\varepsilon}(\mathbf{x}) \coloneqq -k_{\varepsilon} \log(w - P_{\tilde{\Omega}_{\varepsilon}} w)\left(\frac{\mathbf{x}-\mathbf{x}_{\varepsilon}}{k_{\varepsilon}}\right) \quad \text{on } \Omega \; ,
\end{equation}
a fundamental role in the estimates from above is played by the rescaled projection errors $V_{\varepsilon}$ and $\tilde{V}_{\varepsilon}$, defined by the relations:
\begin{equation}\label{VepsilonDefinition}
w = P_{\Omega_{\varepsilon}} w + e^{-\beta_{\varepsilon} \Psi_{\varepsilon}(\mathbf{q})} V_{\varepsilon} \quad \text{on } \Omega_{\varepsilon} \; ,
\end{equation}
\begin{equation}\label{VTildeepsilonDefinition}
w = P_{\tilde{\Omega}_{\varepsilon}} w + e^{-\beta_{\varepsilon} \tilde{\Psi}_{\varepsilon}(\mathbf{x}_{\varepsilon})} \tilde{V}_{\varepsilon} \quad \text{on } \tilde{\Omega}_{\varepsilon} \; .
\end{equation}
We remark that the functions $h_{\varepsilon} \coloneqq e^{-\beta_{\varepsilon} \Psi_{\varepsilon}(\mathbf{q})} V_{\varepsilon}= w - P_{\Omega_{\varepsilon}} w $ and $\tilde{h}_{\varepsilon} \coloneqq e^{-\beta_{\varepsilon} \tilde{\Psi}_{\varepsilon}(\mathbf{x}_{\varepsilon})} \tilde{V}_{\varepsilon}= w - P_{\tilde{\Omega}_{\varepsilon}} w $ solve respectively 
\begin{equation}\label{eqn:ellipticProblemh}
\begin{cases}
-\Delta h_{\varepsilon} + \tilde{\lambda}_0 \underline{m} \, h_{\varepsilon}=0 & \text{in } \Omega_{\varepsilon} \; , \\
h_{\varepsilon} = w & \text{on } \partial\Omega_{\varepsilon} \; ,
\end{cases}
\end{equation}
and
\begin{equation}\label{eqn:ellipticProblemTildeh}
\begin{cases}
-\Delta \tilde{h}_{\varepsilon} + \tilde{\lambda}_0 \underline{m} \, \tilde{h}_{\varepsilon}=0 & \text{in } \tilde{\Omega}_{\varepsilon} \; ,\\
\tilde{h}_{\varepsilon} = w & \text{on } \partial\tilde{\Omega}_{\varepsilon} \; .
\end{cases}
\end{equation}
In particular, by the Strong Maximum Principle, both the functions $V_{\varepsilon}$ and $\tilde{V}_{\varepsilon}$ are strictly positive in their domain of definition.

Moreover, the functions $\Psi_{\varepsilon}$ and $\tilde{\Psi}_{\varepsilon}$ solve respectively
\begin{equation}\label{eqn:psiEqn}
\begin{cases}
- k_{\varepsilon} \Delta \Psi_{\varepsilon} + |\nabla \Psi_{\varepsilon}|^2 - \tilde{\lambda}_0 \underline{m} = 0 & \text{in } \Omega \\
\Psi_{\varepsilon} = -k_{\varepsilon} \log(w)(\frac{\mathbf{x}-\mathbf{q}}{k_{\varepsilon}}) & \text{on }  \partial\Omega
\end{cases}
\end{equation}
and
\begin{equation}\label{eqn:psiTildeEqn}
\begin{cases}
- k_{\varepsilon} \Delta \tilde{\Psi}_{\varepsilon} + |\nabla \tilde{\Psi}_{\varepsilon}|^2 - \tilde{\lambda}_0 \underline{m} = 0 & \text{in } \Omega \\
\tilde{\Psi}_{\varepsilon} = -k_{\varepsilon} \log(w)(\frac{\mathbf{x}-\mathbf{x}_{\varepsilon}}{k_{\varepsilon}}) & \text{on }  \partial\Omega
\end{cases}
\end{equation}
It is clear that a key point in the analysis would be the knowledge of the behavior of $\Psi_{\varepsilon}$, $\tilde{\Psi}_{\varepsilon}$, $V_{\varepsilon}$ and $\tilde{V}_{\varepsilon}$ for $\varepsilon \to 0^{+}$. This is achieved in \cite[Lemma 4.4, Lemma 4.6]{NiWei:Spike95} using vanishing viscosity methods. The particular choice of $f$ influences the proofs only to the extent of the asymptotic behavior of $w(|\mathbf{x}|)$ and $w'(|\mathbf{x}|)$ for $|\mathbf{x}| \to +\infty$, which indeed depends on $f$. This is intuitive if one notices that in \eqref{eqn:ellipticProblemh}, \eqref{eqn:ellipticProblemTildeh}, \eqref{eqn:psiEqn} and \eqref{eqn:psiTildeEqn} there is no explicit dependence on $f$, but only on $w$ through its boundary values, for $|\mathbf{x}| \to +\infty$.

Moreover, comparing Proposition \ref{result:ResultShapeSolutionOverRn} and \cite[Lemma 4.1]{NiWei:Spike95} we can notice that the behavior of $w(|\mathbf{x}|)$ and $w'(|\mathbf{x}|)$ for $|\mathbf{x}| \to +\infty$ are exactly the same, up to a rescaling of $\Omega$.

Thus, eventually rescaling the domain $\Omega$, exactly as in \cite{NiWei:Spike95} we can obtain the following lemmas, of which we omit the proof.
\begin{lemma}[{\cite[Lemma 4.4]{NiWei:Spike95}}]\label{LemmaVepsilon}
(i) $ \Psi_{\varepsilon} \to \Psi_0 \in W^{1, \infty}(\Omega)$ uniformly in $\overline{\Omega}$ as $\varepsilon \to 0$, where $\Psi_0(\mathbf{q}) = 2 \sqrt{\tilde{\lambda}_0 \underline{m}} \, d(\mathbf{q}, \partial\Omega)$ \newline
(ii) For every sequence $\varepsilon_k$ there is a subsequence $\varepsilon_{k_l} \to 0$ such that $V_{\varepsilon_{k_l}} \to V_0$ uniformly on every compact set of $\R^n$, where $V_0$ is a positive solution of
\begin{equation}\label{eq_systemV0}
\begin{cases}
-\Delta u + \tilde{\lambda}_0 \underline{m} \, u = 0 & \text{in } \R^n \; , \\
u(\mathbf{0}) = 1 \; \text{and } u>0  & \text{in } \R^n \; .
\end{cases}
\end{equation}
Moreover, for any $\sigma_1 >0$, $\sup_{\mathbf{x} \in \overline{\Omega}_{\varepsilon_{k_l}}} e^{-(1+\sigma_1) \sqrt{\tilde{\lambda}_0 \underline{m}} |\mathbf{x}|} |V_{\varepsilon_{k_l}}(\mathbf{x}) - V_0(\mathbf{x}) | \to 0$ for $\varepsilon_{k_l} \to 0$.
\end{lemma}
\begin{lemma}[{\cite[Lemma 4.6]{NiWei:Spike95}}]\label{LemmaVTildeepsilon}
(i) There exists a positive constant C such that $ \| \tilde{\Psi}_{\varepsilon} \|_{L^{\infty}(\Omega)} \le C$ \newline
(ii) For any $\sigma_0 >0 $ there is an $\varepsilon_0 > 0$ such that for any $\varepsilon < \varepsilon_0$
\[
\tilde{\Psi}_{\varepsilon}(\mathbf{x}_{\varepsilon}) \le (2+\sigma_0) \sqrt{\tilde{\lambda}_0 \underline{m}} \, d(\mathbf{x}_{\varepsilon}, \partial\Omega)
\]
(iii) For every sequence $\varepsilon_k$ there is a subsequence $\varepsilon_{k_l} \to 0$ such that $\tilde{V}_{\varepsilon_{k_l}} \to \tilde{V}_0$ uniformly on every compact set of $\R^n$, where $\tilde{V}_0$ is a positive solution of
\[
\begin{cases}
-\Delta u + \tilde{\lambda}_0 \underline{m} \, u = 0 \quad \text{in } \R^n \\
u(\mathbf{0}) = 1 \; \text{and } u>0 \; \text{in } \R^n
\end{cases}
\]
Moreover, for any $\sigma_2 >0$, $\sup_{\mathbf{x} \in \overline{\tilde{\Omega}}_{\varepsilon_{k_l}}} e^{-(1+\sigma_2) \sqrt{\tilde{\lambda}_0 \underline{m}} |\mathbf{x}|} |\tilde{V}_{\varepsilon_{k_l}}(\mathbf{x}) - \tilde{V}_0(\mathbf{x}) | \to 0$ for $\varepsilon_{k_l} \to 0$.
\end{lemma}
Notice that Lemma \ref{LemmaVTildeepsilon} (ii) gives an estimate from above for $\tilde{\Psi}_{\varepsilon}(\mathbf{x}_{\varepsilon})$. An estimate from below is contained in the following
\begin{lemma}\label{LemmaVTildeepsilonBoundBelow}
 There exists $ 0 < \sigma_3 < 1$ and $\varepsilon_0 > 0$ such that for any $\varepsilon < \varepsilon_0$
\[
\tilde{\Psi}_{\varepsilon}(\mathbf{x}_{\varepsilon}) \ge \sigma_3 \sqrt{\tilde{\lambda}_0 \underline{m}} \,d(\mathbf{x}_{\varepsilon}, \partial\Omega)
\]
\end{lemma}
\begin{proof}
To begin with, we recall that $\tilde{\Psi}_{\varepsilon}(P_{\varepsilon})$ is a solution of \eqref{eqn:psiTildeEqn}. Moreover, since $\tilde{h}_{\varepsilon} \coloneqq e^{-\beta_{\varepsilon} \tilde{\Psi}_{\varepsilon}(\mathbf{x}_{\varepsilon})} \tilde{V}_{\varepsilon}= w - P_{\tilde{\Omega}_{\varepsilon}} w $ solves \eqref{eqn:ellipticProblemTildeh} and
\begin{equation}\label{wAsymptExpansion}
w(|\mathbf{x}|) \sim \frac{C}{|\mathbf{x}|} e^{-\sqrt{\tilde{\lambda}_0 \underline{m}} \, |\mathbf{x}|} \quad \text{for } |\mathbf{x}| \to +\infty, \, C>0
\end{equation}
thanks to Proposition \ref{result:ResultShapeSolutionOverRn}, after standard computations it also holds that (see also Lemma 4.5 in \cite{NiWei:Spike95})
\begin{equation}\label{BoundaryPsiTilveEpsilon}
\lim_{\varepsilon \to 0} \frac{\tilde{\Psi}_{\varepsilon}(\mathbf{x})}{|\mathbf{x}-\mathbf{x}_{\varepsilon}|} = \sqrt{\tilde{\lambda}_0 \underline{m}},
\qquad\text{uniformly for }\mathbf{x} \in \partial\Omega.
\end{equation}
Let us introduce a point $\bar{\mathbf{x}}_{\varepsilon} \in \partial\Omega$ such that $d(\mathbf{x}_{\varepsilon},\partial\Omega) = |\mathbf{x}_{\varepsilon} - \bar{\mathbf{x}}_{\varepsilon}|$ and a point $\mathbf{y}_{\varepsilon} \in \Omega^c$ that will be positioned sufficiently close to the boundary on the line passing through $\mathbf{x}_{\varepsilon}$ and $\bar{\mathbf{x}}_{\varepsilon}$. Then, we can define
\[
v_{\varepsilon}(\mathbf{x}) \coloneqq \sigma_3 \sqrt{\tilde{\lambda}_0 \underline{m}} \,|\mathbf{x}-\mathbf{y}_{\varepsilon}| \in C^2(\overline{\Omega}) \; .
\]
Observe that
\[
v_{\varepsilon}(\mathbf{x}_{\varepsilon}) = \sigma_3 \sqrt{\tilde{\lambda}_0 \underline{m}} \,|\mathbf{x}_{\varepsilon}-\mathbf{y}_{\varepsilon}| = \sigma_3 \sqrt{\tilde{\lambda}_0 \underline{m}} \,(|\mathbf{x}_{\varepsilon} - \bar{\mathbf{x}}_{\varepsilon}| + |\bar{\mathbf{x}}_{\varepsilon} - \mathbf{y}_{\varepsilon}|) \ge
\]
\begin{equation}\label{EstimateV1}
\ge \sigma_3 \sqrt{\tilde{\lambda}_0 \underline{m}} \,|\mathbf{x}_{\varepsilon} - \bar{\mathbf{x}}_{\varepsilon}| = \sigma_3 \sqrt{\tilde{\lambda}_0 \underline{m}} \, d(\mathbf{x}_{\varepsilon}, \partial\Omega) ; .
\end{equation}
Moreover, by direct computation:
\[
\nabla v_{\varepsilon}(\mathbf{x}) = \sigma_3 \sqrt{\tilde{\lambda}_0 \underline{m}} \frac{\mathbf{x}-\mathbf{y}_{\varepsilon}}{|\mathbf{x}-\mathbf{y}_{\varepsilon}|}
\]
and
\[
\Delta v_{\varepsilon}(\mathbf{x})= \sigma_3 \sqrt{\tilde{\lambda}_0 \underline{m}} \frac{(N-1)}{|\mathbf{x}-\mathbf{y}_{\varepsilon}|},
\]
where $N$ is the dimension of the space. Hence:
\begin{equation}\label{EstimateV2}
- k_{\varepsilon} \Delta v_{\varepsilon}(\mathbf{x}) + |\nabla v_{\varepsilon}(\mathbf{x})|^2 - \tilde{\lambda}_0 \underline{m} =- k_{\varepsilon} \sigma_3 \sqrt{\tilde{\lambda}_0 \underline{m}} \frac{(N-1)}{|\mathbf{x}-\mathbf{y}_{\varepsilon}|} +  \sigma_3^2 \tilde{\lambda}_0 \underline{m} - \tilde{\lambda}_0 \underline{m} < 0
\end{equation}
if $\sigma_3$ is chosen sufficiently small. Now, since
\[
v_{\varepsilon}(\mathbf{x}) \le \sigma_3 \sqrt{\tilde{\lambda}_0 \underline{m}} (|\mathbf{x}-\mathbf{x}_{\varepsilon}| + |\mathbf{x}_{\varepsilon} - \bar{\mathbf{x}}_{\varepsilon}| + |\bar{\mathbf{x}}_{\varepsilon} - \mathbf{y}_{\varepsilon}|)
\]
and taking into account \eqref{BoundaryPsiTilveEpsilon}, if
\begin{equation}\label{ConditionV1}
\sigma_3 (|\mathbf{x}_{\varepsilon} - \bar{\mathbf{x}}_{\varepsilon}| + |\bar{\mathbf{x}}_{\varepsilon} - \mathbf{y}_{\varepsilon}|) \le (1+o(1) - \sigma_3) |\mathbf{x}-\mathbf{x}_{\varepsilon}| \; ,
\end{equation}
we have that 
\begin{equation}\label{EstimateV3}
\tilde{\Psi}_{\varepsilon} \ge v_{\varepsilon} \; \text{on } \partial\Omega \; .
\end{equation} 
Since $|\mathbf{x}_{\varepsilon} - \bar{\mathbf{x}}_{\varepsilon}| \le |\mathbf{x}-\mathbf{x}_{\varepsilon}|$ for $\mathbf{x} \in \partial\Omega$, if we choose $\mathbf{y}_{\varepsilon}$ such that $|\bar{\mathbf{x}}_{\varepsilon} - \mathbf{y}_{\varepsilon}| \le |\mathbf{x}_{\varepsilon} - \bar{\mathbf{x}}_{\varepsilon}|$ condition \eqref{ConditionV1} is satisfied for $0 < \sigma_3 < 1/3$ for $\varepsilon$ sufficiently small, and condition \eqref{EstimateV3} is satisfied too.

Combining \eqref{EstimateV2}, \eqref{EstimateV3} and using the Comparison Principle for Quasilinear Equations (see e.g. \cite[Paragraph 10.1]{GilbargTrudinger:Elliptic98}) gives $\tilde{\Psi}_{\varepsilon} \ge v_{\varepsilon} $ on $\Omega$. Relation \eqref{EstimateV1} then concludes the proof.
\end{proof}
To obtain the estimates from above, the following lemma will be useful:
\begin{lemma}\label{LemmaBoundsErrorEstimateTerms}
There exists $D>0$ independent on $\sigma_1$ such that, for any sequence $\varepsilon_n$ there is a subsequence $\varepsilon_{n_k}$ such that:
\begin{equation} \label{EstimateTerm1}
\int_{\Omega_{\varepsilon_{n_k}}} |2w - e^{-\beta_{\varepsilon_{n_k}} \Psi_{\varepsilon_{n_k}}(\mathbf{q})}  V_{\varepsilon_{n_k}}|  |V_{\varepsilon_{n_k}}| = o\left(e^{D \sigma_1 \beta_{\varepsilon_{n_k}}\Psi_{\varepsilon_{n_k}}(\mathbf{q})}\right)
\end{equation}
Moreover, there exists $C>0$ such that:
\begin{equation} \label{EstimateTerm2}
0\le \int_{\Omega_{\varepsilon_{n_k}}} f_0(w) V_{\varepsilon_{n_k}} \le C
\end{equation}
\begin{equation} \label{EstimateTerm3}
0\le \int_{\Omega_{\varepsilon_{n_k}}} f_0( V_{\varepsilon_{n_k}}) V_{\varepsilon_{n_k}} \le C.
\end{equation}
\end{lemma}
\begin{proof}
For the sake of readibility, we do not write $\varepsilon$ with the subscripts. However, we remark that we need to pass to subsequences to apply Lemma \ref{LemmaVepsilon} (ii) and Lemma \ref{LemmaVTildeepsilon} (iii). This remark holds also for all the subsequent proofs, and will be omitted.
 
We recall preliminarly that, since $P_{\Omega_{\varepsilon}} w =  w - e^{-\beta_{\eps} \Psi_{\varepsilon}(\mathbf{q})} V_{\varepsilon}$ solves \eqref{eq:Pw}, it is 
non-negative on $\Omega_{\varepsilon}$ by the Maximum Principle. Hence it holds that
\begin{equation}\label{bound1}
 |2w - e^{-\beta_{\varepsilon} \Psi_{\varepsilon}(\mathbf{q})}  V_{\varepsilon} | = w+P_{\Omega_{\varepsilon}} w  \le 2 |w|.
\end{equation}
 Now, to prove \eqref{EstimateTerm1} we notice that from \eqref{wAsymptExpansion}, \eqref{bound1} and Lemma \ref{LemmaVepsilon}
\begin{align*}
&\int_{\Omega_{\varepsilon}}|2w - e^{-\beta_{\varepsilon} \Psi_{\varepsilon}(\mathbf{q})}  V_{\varepsilon}|  |V_{\varepsilon}|
 \le C \int_{\Omega_{\varepsilon}}|w|  e^{(1+\sigma_1) \sqrt{\tilde{\lambda}_0 \underline{m}}|\mathbf{x}|} \le \\
& \le  C e^{2 b \, \sigma_1 \beta_{\varepsilon}  \sqrt{\tilde{\lambda}_0 \underline{m}} \,  d(\mathbf{q}, \partial\Omega)} = C e^{(b+o(1)) \, \sigma_1 \beta_{\varepsilon} \Psi_{\varepsilon}(\mathbf{q})} 
\end{align*}
with $b \coloneqq \diam(\Omega)/d(\mathbf{q}, \partial\Omega)$, which gives \eqref{EstimateTerm1} for $D = 3b > 2b+o(1)$ for $\varepsilon$ sufficiently small.

Relations \eqref{EstimateTerm2} and \eqref{EstimateTerm3} follow from the fact that $f_0(w)$ and $f_0(V_{\varepsilon})$ have fixed compact support in $\R^n$ with respect to $\varepsilon$, from Lemma \ref{LemmaVepsilon} and the fact that  by elliptic regularity $w \in C^{1}(\R^n)$ and $V_0 \in C^{\infty}(\R^n)$.
\end{proof}
We now proceed to prove the following estimate from above:
\begin{proposition}\label{AboveErrorEstimate}
For any sequence $\varepsilon_n \to 0$ there is a subsequence $\varepsilon_{n_k} \to 0$ such that:
\begin{equation}\label{EstimateFromAbove}
\tilde{\lambda}_{\varepsilon_{n_k}} \le \tilde{\lambda}_0 + \Phi \, e^{-\beta_{\varepsilon_{n_k}} \Psi_{\varepsilon_{n_k}}(\mathbf{q})} + o\biggl(e^{-\beta_{\varepsilon_{n_k}} \Psi_{\varepsilon_{n_k}}(\mathbf{q})}\biggr) \quad \text{for } k \to +\infty \; ,
\end{equation}
where 
\[
\Phi \coloneqq \frac{2 \gamma}{\int_{\R^n} \tilde{m}_0 \, w^2}
\]
and $\gamma \coloneqq \int_{\R^n} f_0(w) \, V_0 > 0$ is independent on the particular $V_0$ solution of \eqref{eq_systemV0}.
\end{proposition}
\begin{proof}
The estimate \eqref{EstimateFromAbove} is obtained from the inequality
\[
\tilde{\lambda}_{\varepsilon} \le \frac{\int_{\Omega_{\varepsilon}} | \nabla P_{\Omega_{\varepsilon}}w|^2}{\int_{\Omega_{\varepsilon}} \tilde{m}_0 \, |P_{\Omega_{\varepsilon}}w|^2 }
\]
Using \eqref{eq:Pw} and \eqref{VepsilonDefinition} we can write:
\[
\begin{split}
\int_{\Omega_{\varepsilon}} | \nabla P_{\Omega_{\varepsilon}}w|^2 &= \int_{\Omega_{\varepsilon}}f_0(w) \, P_{\Omega_{\varepsilon}}w - \int_{\Omega_{\varepsilon}}\tilde{\lambda}_0 \underline{m} |P_{\Omega_{\varepsilon}} w|^2 = \\
& = \int_{\Omega_{\varepsilon}}f_0(w) \, (w - e^{-\beta_{\varepsilon} \Psi_{\varepsilon}(\mathbf{q})} V_{\varepsilon})  - \int_{\Omega_{\varepsilon}}\tilde{\lambda}_0 \underline{m} (w^2-2 e^{-\beta_{\varepsilon} \Psi_{\varepsilon}(\mathbf{q})} w V_{\varepsilon} + e^{-2 \beta_{\varepsilon} \Psi_{\varepsilon}(\mathbf{q})} V_{\varepsilon}^2) = \\
& = \tilde{\lambda}_0 \int_{\Omega_{\varepsilon}}\tilde{m}_0\, w^2 -  e^{-\beta_{\varepsilon} \Psi_{\varepsilon}(\mathbf{q})}\int_{\Omega_{\varepsilon}} f_0(w) V_{\varepsilon} +  \tilde{\lambda}_0 \underline{m} e^{-\beta_{\varepsilon} \Psi_{\varepsilon}(\mathbf{q})} \int_{\Omega_{\varepsilon}} (2w - e^{-\beta_{\varepsilon} \Psi_{\varepsilon}(\mathbf{q})}  V_{\varepsilon})  V_{\varepsilon} 
\end{split}
\]
and 
\begin{align*}\label{eq:m0Pw}
& \frac{1}{\int_{\Omega_{\varepsilon}} \tilde{m}_0 \, |P_{\Omega_{\varepsilon}}w|^2} = 
\frac{1}{\int_{\Omega_{\varepsilon}} \tilde{m}_0 \, (w^2-2 e^{-\beta_{\varepsilon} \Psi_{\varepsilon}(\mathbf{q})} w V_{\varepsilon} + e^{-2 \beta_{\varepsilon} \Psi_{\varepsilon}(\mathbf{q})} V_{\varepsilon}^2)} = \\
& = \frac{1}{\int_{\Omega_{\varepsilon}} \tilde{m}_0 \, w^2} \left[1 +\frac{e^{-\beta_{\varepsilon} \Psi_{\varepsilon}(\mathbf{q})}}{\int_{\Omega_{\varepsilon}} \tilde{m}_0 \, w^2}  \int_{\Omega_{\varepsilon}} \tilde{m}_0 \left(2w - e^{-\beta_{\varepsilon} \Psi_{\varepsilon}(\mathbf{q})}  V_{\varepsilon}\right)  V_{\varepsilon} + o\left(e^{-\beta_{\varepsilon} \Psi_{\varepsilon}(\mathbf{q})}\right) \right],
\end{align*}
where in the last passage we have used \eqref{EstimateTerm1}, with $\sigma_1$ small enough. Combining these results with \eqref{EstimateTerm1} and \eqref{EstimateTerm2} we obtain the following estimate from above for $\tilde{\lambda}_{\varepsilon}$:
\[
\tilde{\lambda}_{\varepsilon} \le \tilde{\lambda}_0 + \frac{1}{\int_{\Omega_{\varepsilon}} \tilde{m}_0 \, w^2} \biggl[ e^{-\beta_{\varepsilon} \Psi_{\varepsilon}(\mathbf{q})}\int_{\Omega_{\varepsilon}} f_0(w) V_{\varepsilon}  - e^{-2 \beta_{\varepsilon} \Psi_{\varepsilon}(\mathbf{q})}\int_{\Omega_{\varepsilon}} f_0(V_{\varepsilon}) V_{\varepsilon} + o (e^{-\beta_{\varepsilon} \Psi_{\varepsilon}(\mathbf{q})})\biggr].
\]
Now, using \eqref{EstimateTerm3} and the facts that 
\[
\int_{\Omega_{\varepsilon}} \tilde{m}_0 \, w^2 \to \int_{\R^n} \tilde{m}_0 \, w^2 > 0 \quad \text{for } \varepsilon \to 0
\]
\[
\int_{\Omega_{\varepsilon}} f_0(w) V_{\varepsilon} \to \int_{\R^n} f_0(w) V_0 \quad \text{for } \varepsilon \to 0
\]
from the previous estimate from above we obtain
\[
\tilde{\lambda}_{\varepsilon} \le \tilde{\lambda}_0 + \frac{e^{-\beta_{\varepsilon} \Psi_{\varepsilon}(\mathbf{q})}}{\int_{\R^n} \tilde{m}_0 \, w^2} \int_{\R^n} f_0(w) V_0 + o (e^{-\beta_{\varepsilon} \Psi_{\varepsilon}(\mathbf{q})})
\]
which is \eqref{EstimateFromAbove}. The fact that $\gamma$ is positive and independent on the particular $V_0$ can be proved exactly as in \cite{NiWei:Spike95} Lemma 4.7. 
\end{proof}

In view of the estimate from below, we also prove an estimate from above depending on $\tilde{\Psi}_{\varepsilon}(\mathbf{x}_{\varepsilon})$. 
\begin{proposition}\label{AboveErrorEstimate2}
For any sequence $\varepsilon_n \to 0$ there is a subsequence $\varepsilon_{n_k} \to 0$ such that:
\begin{equation}\label{EstimateFromAbove2}
\tilde{\lambda}_{\varepsilon_{n_k}} \le \tilde{\lambda}_0 + \Phi \, e^{-\beta_{\varepsilon_{n_k}} \tilde{\Psi}_{\varepsilon_{n_k}}(\mathbf{x}_{\varepsilon_{n_k}})} + o\biggl(e^{-\beta_{\varepsilon_{n_k}} \tilde{\Psi}_{\varepsilon_{n_k}}(\mathbf{x}_{\varepsilon_{n_k}})} \biggr) \quad \text{for } k \to +\infty \; ,
\end{equation}
where 
\[
\Phi \coloneqq \frac{2 \gamma}{\int_{\R^n} \tilde{m}_0 \, w^2}
\]
and $\gamma \coloneqq \int_{\R^n} f_0(w) \, \tilde{V}_0 > 0$ is independent on the particular $\tilde{V}_0$ solution of \eqref{eq_systemV0}.
\end{proposition}
\begin{proof}
Similarly as in the proof of Proposition \ref{AboveErrorEstimate}, we have:
\begin{equation}\label{LembdaEpsilonAboveErrorEstimate2}
\tilde{\lambda}_{\varepsilon} \le \frac{\int_{\tilde{\Omega}_{\varepsilon}} | \nabla P_{\tilde{\Omega}_{\varepsilon}}w|^2}{\int_{\tilde{\Omega}_{\varepsilon}} \tilde{m}_0 \, |P_{\tilde{\Omega}_{\varepsilon}}w|^2 } = \frac{\tilde{\lambda}_0 + B}{1-A}
\end{equation}
where
\[
\begin{split}
A &:=  \frac{ e^{-\beta_{\varepsilon} \tilde{\Psi}_{\varepsilon}(\mathbf{x}_{\varepsilon})}}{\int_{\tilde{\Omega}_{\varepsilon}} \tilde{m}_0 \, w^2} \int_{\tilde{\Omega}_{\varepsilon}} \tilde{m}_0 (2w - e^{-\beta_{\varepsilon} \tilde{\Psi}_{\varepsilon}(\mathbf{x}_{\varepsilon})}  \tilde{V}_{\varepsilon})  \tilde{V}_{\varepsilon}  \\
B &:=  \frac{1}{\int_{\tilde{\Omega}_{\varepsilon}} \tilde{m}_0 \, w^2}  \biggl[\tilde{\lambda}_0 \underline{m} e^{-\beta_{\varepsilon} \tilde{\Psi}_{\varepsilon}(\mathbf{x}_{\varepsilon})} \int_{\tilde{\Omega}_{\varepsilon}} (2w - e^{-\beta_{\varepsilon} \tilde{\Psi}_{\varepsilon}(\mathbf{x}_{\varepsilon})}  \tilde{V}_{\varepsilon})  \tilde{V}_{\varepsilon}   -  e^{-\beta_{\varepsilon} \tilde{\Psi}_{\varepsilon}(\mathbf{x}_{\varepsilon})}\int_{\tilde{\Omega}_{\varepsilon}} f_0(w) \tilde{V}_{\varepsilon}   \biggr]
\end{split}
\]
so that the sum $(\tilde{\lambda}_0 A+B)$ is of order $ e^{-\beta_{\varepsilon} \tilde{\Psi}_{\varepsilon}(\mathbf{x}_{\varepsilon})}$, since it satisfies
\[
(\tilde{\lambda}_0 A+B) = \frac{ e^{-\beta_{\varepsilon} \tilde{\Psi}_{\varepsilon}(\mathbf{x}_{\varepsilon})}}{\int_{\tilde{\Omega}_{\varepsilon}} \tilde{m}_0 \, w^2} \int_{\tilde{\Omega}_{\varepsilon}} f_0(w) \tilde{V}_{\varepsilon}  + \frac{ e^{-2\beta \tilde{\Psi}_{\varepsilon}(\mathbf{x}_{\varepsilon})}}{\int_{\tilde{\Omega}_{\varepsilon}} \tilde{m}_0 \, w^2} \int_{\tilde{\Omega}_{\varepsilon}} f_0(\tilde{V}_{\varepsilon}) \tilde{V}_{\varepsilon}  
\]
with
\[
\int_{\tilde{\Omega}_{\varepsilon}} \tilde{m}_0 \, w^2 \to \int_{\R^n} \tilde{m}_0 \, w^2 \qquad \text{for } \varepsilon \to 0
\]
\[
\int_{\tilde{\Omega}_{\varepsilon}} f_0(w) \tilde{V}_{\varepsilon} \to \int_{\R^n} f_0(w) \tilde{V}_0 \qquad \text{for } \varepsilon \to 0
\]
\[
 \int_{\tilde{\Omega}_{\varepsilon}} f_0(\tilde{V}_{\varepsilon}) \tilde{V}_{\varepsilon}  \to  \int_{\R^n} f_0(\tilde{V}_0) \tilde{V}_0  \qquad \text{for } \varepsilon \to 0.
\]
This implies
\begin{equation}\label{DominantTermAboveErrorEstimate2}
(\tilde{\lambda}_0 A+B) = \Phi e^{-\beta_{\eps} \tilde{\Psi}_{\varepsilon}(\mathbf{x}_{\varepsilon})} + o(e^{-\beta_{\eps} \tilde{\Psi}_{\varepsilon}(\mathbf{x}_{\varepsilon})})  \qquad \text{for } \varepsilon \to 0.
\end{equation}
Moreover, notice that  $A = o(1)$ and $B = o(1)$ for $\varepsilon \to 0$. Now, since
\begin{equation}\label{DominantTermAboveErrorEstimate3}
\frac{\tilde{\lambda}_0 + B}{1-A} = \tilde{\lambda}_0 + (\tilde{\lambda}_0 A + B) \frac{1}{1-A} =  \tilde{\lambda}_0 + (\tilde{\lambda}_0 A + B) \sum_{k=0}^{+\infty} A^k
\end{equation}
the first order term in the expansion of $\frac{\int_{\tilde{\Omega}_{\varepsilon}} | \nabla P_{\tilde{\Omega}_{\varepsilon}}w|^2}{\int_{\tilde{\Omega}_{\varepsilon}} \tilde{m}_0 \, |P_{\tilde{\Omega}_{\varepsilon}}w|^2 } = \frac{\tilde{\lambda}_0 + B}{1-A}$ is exactly $\tilde{\lambda}_0 A + B$. Hence, combining \eqref{LembdaEpsilonAboveErrorEstimate2}, \eqref{DominantTermAboveErrorEstimate3} and \eqref{DominantTermAboveErrorEstimate2} we conclude.
\end{proof}
\begin{remark}\label{rem:proofabove}
Notice that, as a consequence of the proof above, we have
\[
\frac{\int_{\tilde{\Omega}_{\varepsilon}} | \nabla  P_{\tilde{\Omega}_{\varepsilon}}|^2}{\int_{\tilde{\Omega}_{\varepsilon}} \tilde{m}_0 \, | P_{\tilde{\Omega}_{\varepsilon}}|^2 } = \tilde{\lambda}_0 + \Phi e^{-\beta_{\varepsilon} \tilde{\Psi}_{\varepsilon}(\mathbf{x}_{\varepsilon})} + o(e^{-\beta_{\varepsilon} \tilde{\Psi}_{\varepsilon}(\mathbf{x}_{\varepsilon})}).
\]
\end{remark}

\section{Estimate from below}\label{section:EstimateBelow}

To prove the bound from below, it is convenient to choose a different normalization 
of the blow-up sequence. More precisely, with a little abuse of notation, throughout this section we denote with $\tilde u_{ \varepsilon}$ the family of solutions to the blow up problem, normalized in such  a way that $\tilde u_{\varepsilon}(\mathbf{0})=1$.
For any sequence $\varepsilon_n \to 0$, this family has the same qualitative properties of the one with $L^2(\tilde{\Omega}_{\varepsilon})$ normalization; in particular, up to a further subsequence, Proposition \ref{prop:PropositionBlowUp} still holds if we consider $w$ normalized such that $w(\mathbf{0}) = 1$. 

Let now $\tilde u_{ \varepsilon} = P_{\tilde{\Omega}_{\varepsilon}} w + e^{-\beta_{\varepsilon} \tilde{\Psi}_{\varepsilon}(\mathbf{x}_{\varepsilon})} \phi_{\varepsilon}$, so that $\phi_{\varepsilon}$ is defined by:
\begin{equation}\label{eq:phi}
\phi_{\varepsilon} \coloneqq e^{\beta_{\varepsilon} \tilde{\Psi}_{\varepsilon}(\mathbf{x}_{\varepsilon})} (\tilde u_{\varepsilon} - P_{\tilde{\Omega}_{\varepsilon}} w)
\end{equation}
Our aim is to prove the following
\begin{proposition}\label{PhiErrorEstimate}
For every sequence $\varepsilon_n \to 0$ there exists a subsequence $\varepsilon_{n_k}$ such that
\[ 
\| \phi_{\varepsilon_{n_k}} \|_{H^2(\tilde{\Omega}_{\varepsilon_{n_k}})} \le C. 
\] 
\end{proposition}
Before giving the proof we need some preparation. We remark that, if we introduce $\hat{u}_{ \varepsilon}$ as the blow up sequence with normalizion $\hat{u}_{ \varepsilon}(\mathbf{0}) = P_{\tilde{\Omega}_{\varepsilon}} w (\mathbf{0})$, from Lemma \ref{LemmaVTildeepsilon} and the definition of $\tilde{V}_{\varepsilon}$ we have that, up to a subsequence $\tilde u_{ \varepsilon} - \hat{u}_{ \varepsilon} = C_{\varepsilon} \, e^{-\beta_{\eps} \tilde{\Psi}_{\varepsilon}(\mathbf{x}_{\varepsilon})} \tilde u_{ \varepsilon}$ (of course we are considering $w$ s.t. $w(\mathbf{0})=1$), with $|C_{\varepsilon}| \le C$. This in turn means that:
\begin{align*}
\phi_{\varepsilon} &= e^{\beta_{\varepsilon} \tilde{\Psi}_{\varepsilon}(\mathbf{x}_{\varepsilon})} (\tilde u_{\varepsilon} - P_{\tilde{\Omega}_{\varepsilon}} w) = e^{\beta_{\varepsilon} \tilde{\Psi}_{\varepsilon}(\mathbf{x}_{\varepsilon})} (\tilde u_{\varepsilon} - \hat{u}_{ \varepsilon} + \hat{u}_{ \varepsilon} - P_{\tilde{\Omega}_{\varepsilon}} w) = \\
& = C_{\varepsilon} \, \tilde  u_{ \varepsilon} + e^{\beta_{\eps} \tilde{\Psi}_{\varepsilon}(\mathbf{x}_{\varepsilon})} (\hat{u}_{\varepsilon} - P_{\tilde{\Omega}_{\varepsilon}} w)
\end{align*}
and if we define $\hat{\phi}_{\varepsilon} \coloneqq e^{\beta_{\eps} \tilde{\Psi}_{\varepsilon}(\mathbf{x}_{\varepsilon})} (\hat{u}_{\varepsilon} - P_{\tilde{\Omega}_{\varepsilon}} w)$ we have
\begin{equation}\label{PhiEstimateWRTPhiTilde}
\phi_{\varepsilon} = C_{\varepsilon} \,\tilde  u_{ \varepsilon} + \hat{\phi}_{\varepsilon}
\end{equation}
so that what we actually need to estimate is $\hat{\phi}_{\varepsilon}$. Moreover, observe that by definition
\begin{equation}\label{PhiTildeNull}
\hat{\phi}_{\varepsilon}(\mathbf{0}) = 0
\end{equation}
and this will be a key property in the following. 

If we define $\tilde{\lambda}_1 \coloneqq \tilde{\lambda}_{\varepsilon} - \tilde{\lambda}_0$, it is easily seen that the function $\hat{\phi}_{\varepsilon}$ solves the following problem:
\begin{equation}
\begin{cases}\label{eq:PhiHat}
-\Delta \hat{\phi}_{\varepsilon} + \tilde{\lambda}_0 \underline{m} \hat{\phi}_{\varepsilon} = f_0(\hat{\phi}_{\varepsilon}) - f_0(\tilde{V}_{\varepsilon}) + e^{\beta_{\varepsilon} \tilde{\Psi}_{\varepsilon}(\mathbf{x}_{\varepsilon})} \tilde{\lambda}_1 \tilde{m}_0 \hat{u}_{\varepsilon} & \text{in } \tilde{\Omega}_{\varepsilon} \; , \\
\hat{\phi}_{\varepsilon} = 0 & \text{on } \partial\tilde{\Omega}_{\varepsilon} \; .
\end{cases}
\end{equation}
Moreover, notice that from Proposition \ref{AboveErrorEstimate2} we know that $ | e^{\beta_{\varepsilon} \tilde{\Psi}_{\varepsilon}(\mathbf{x}_{\varepsilon})} \tilde{\lambda}_1| \le C $.

Following the ideas in \cite{NiWei:Spike95} we first prove the following lemma:
\begin{lemma}\label{PhiHatRHS}
For every sequence $\varepsilon_n \to 0$ there exists a subsequence $\varepsilon_{n_k}$ 
and a solution $\tilde{V}_0$ of \eqref{eq_systemV0} such that
\begin{align*}
 \| f_0(\tilde{V}_0) - f_0(\tilde{V}_{\varepsilon_{n_k}}) + e^{\beta_{\varepsilon_{n_k}} \tilde{\Psi}_{\varepsilon_{n_k}}(\mathbf{x}_{\varepsilon_{n_k}})} \tilde{\lambda}_1 \tilde{m}_0 (\hat{u}_{\varepsilon_{n_k}} -w) \|_{L^2(\tilde{\Omega}_{\varepsilon_{n_k}})} \le o(1)
\end{align*}
for $k \to +\infty$.
\end{lemma}
\begin{proof}
Since, by Lemma \ref{LemmaVTildeepsilon} (iii), up to a subsequence $\tilde{V}_{\varepsilon} \to \tilde{V}_0$ uniformly on compact sets of $\R^n$ and $\| \hat{u}_{\varepsilon} -w \|_{L^2(\R^n)} \to 0$ we have that:
\begin{align*}
 \| f_0(\tilde{V}_0) - f_0(\tilde{V}_{\varepsilon}) + e^{\beta_{\varepsilon} \tilde{\Psi}_{\varepsilon}(\mathbf{x}_{\varepsilon})} \tilde{\lambda}_1 \tilde{m}_0 (\hat{u}_{\varepsilon} -w) \|_{L^2(\tilde{\Omega}_{\varepsilon})} \le \\
 \le \| f_0(\tilde{V}_0) - f_0(\tilde{V}_{\varepsilon}) \|_{L^2(\tilde{\Omega}_{\varepsilon})} + C \| \tilde{m}_0 (\hat{u}_{\varepsilon} -w)\|_{L^2(\tilde{\Omega}_{\varepsilon})}  \le  o(1).
\end{align*}
\end{proof}
Proposition \ref{PhiErrorEstimate} will follow from the following estimate for $\hat{\phi}_{\varepsilon}$:
\begin{lemma}\label{Lemma:PhiHatW2Estimate}
For every sequence $\varepsilon_n \to 0$ there exists a subsequence $\varepsilon_{n_k}$ such that
\[ 
\| \hat{\phi}_{\varepsilon_{n_k}} \|_{H^2(\tilde{\Omega}_{\varepsilon})} \le C
\] 
\end{lemma}
\begin{proof}
We follow the proof of \cite{NiWei:Spike95} Lemma 6.7, Corollary 6.8. In particular we begin proving by contradiction that $\| \hat{\phi}_{\varepsilon} \|_{L^{2}(\tilde{\Omega}_{\varepsilon})} \le C$. 

Suppose that there is a sequence such that $\| \hat{\phi}_{\varepsilon_k} \|_{L^{2}(\tilde{\Omega}_{\varepsilon_k})} \to +\infty$ for $\varepsilon_k \to 0$. We define:
\[ M_k \coloneqq \| \hat{\phi}_{\varepsilon_k} \|_{L^{2}(\tilde{\Omega}_{\varepsilon_k})} ,\quad g_{k} \coloneqq \frac{\hat{\phi}_{\varepsilon_k}}{M_k}, \quad
F_k \coloneqq   e^{\beta_{\varepsilon_k} \tilde{\Psi}_{\varepsilon_k}(\mathbf{x}_{\varepsilon_k})} \tilde{\lambda}_1 \tilde{m}_0 \hat{u}_{\varepsilon_k} -  f_0(\tilde{V}_{\varepsilon_k}) \]
so that $\| g_k \|_{L^{2}(\tilde{\Omega}_{\varepsilon_k})} = 1$. We observe that $g_k$ solves
\[
\begin{cases}
-\Delta g_k + \tilde{\lambda}_0 \underline{m} g_k= f_0(g_k) + \frac{1}{M_k} F_k \quad \text{in } \tilde{\Omega}_{\varepsilon_k} \\
g_k = 0 \quad \text{on } \partial\tilde{\Omega}_{\varepsilon_k}
\end{cases}
\]
Hence, from \cite[Lemma 6.4]{NiWei:Spike95} (with $\mu=0$ and $s=2$) we have that:
\[
\| g_k \|_{H^2(\tilde{\Omega}_{\varepsilon_k})} \le C \bigl( \| f_0(g_k) \|_{L^{2}(\tilde{\Omega}_{\varepsilon_k})} + \| F_k/M_k \|_{L^{2}(\tilde{\Omega}_{\varepsilon_k})} \bigl), 
\]
where $C$ is independent of $0<\eps\le\eps_0$. 
Now, from Lemma \ref{PhiHatRHS} we know that:
\begin{align*}
& \| F_k/M_k \|_{L^{2}(\tilde{\Omega}_{\varepsilon_k})} =
\frac{1}{M_k} \| f_0(\tilde{V}_0) - e^{\beta_{\varepsilon} \tilde{\Psi}_{\varepsilon_k}(\mathbf{x}_{\varepsilon_k})} \tilde{\lambda}_1 \tilde{m}_0 w \|_{L^{2}(\tilde{\Omega}_{\varepsilon_k})} \le \\
& + \frac{1}{M_k} \| f_0(\tilde{V}_0) \|_{L^{2}(\tilde{\Omega}_{\varepsilon_k})} + \frac{C}{M_k} \| w \|_{L^{2}(\tilde{\Omega}_{\varepsilon_k})} = o(1) \; ,
\end{align*}
since $\tilde{V}_0$ is bounded on bounded sets of $\R^n$ and $ \| w \|_{L^{2}(\R^n)}$ is finite. Moreover
\[
 \| f_0(g_k) \|_{L^{2}(\tilde{\Omega}_{\varepsilon_k})} \le C \; ,
\]
from the normalization chosen for $g_k$. 

From the previous discussion we deduce that:
\begin{equation}\label{g_kBoundedW2Norm}
\| g_k \|_{H^2(\tilde{\Omega}_{\varepsilon_k})} \le C \bigl( \| f_0(g_k) \|_{L^{2}(\tilde{\Omega}_{\varepsilon_k})} + o(1) \bigr) \le C \; .
\end{equation}
This means that we can extract a further subsequence (still called $g_k$) that converges weakly in $H^1(\R^N)$ and strongly in $C^{1, \alpha}_{loc}(\R^N)$ for some $0 < \alpha < 1$, to a function $g \in H^1(\R^N)$ solution of:
\[
-\Delta g + \tilde{\lambda}_0 \underline{m} \, g= f_0(g)  \quad \text{in } \R^n.
\]
From Proposition \ref{result:ResultProblemOverRn}, we know that the only 
$H^1(\R^N)$-solutions to such equation are multiples of $w$, so that:
\[
g = a \, w \qquad \text{for  some }a \in \R.
\]
However, we know that $w(\mathbf{0})>0$ while from \eqref{PhiTildeNull} we have that $g_k(\mathbf{0}) = \hat{\phi}_{\varepsilon_k}(\mathbf{0})/M_k = 0$, so that $a = 0$ and:
\begin{equation}\label{PhiHatVanishingWeakLimit}
g \equiv 0.
\end{equation}
But now we state that $\| g_k \|_{H^2(\tilde{\Omega}_{\varepsilon_k})} \to 0$ for $\varepsilon_k \to 0$, reaching a contradiction since by choice  $\| g_k \|_{L^{2}(\tilde{\Omega}_{\varepsilon_k})} = 1$. 

Indeed, observe that from the strong convergence of $g_k$ to $g=0$ on compact subsets of $\R^n$ we have that
\[
\| f_0(g_k) \|_{L^{2}(\tilde{\Omega}_{\varepsilon_k})} = \| f_0(g_k) \|_{L^{2}(B^1(\mathbf{0}))} = o(1) \; \text{for } \varepsilon_k \to 0
\]
and combining this with \eqref{g_kBoundedW2Norm} gives $\| g_k \|_{H^2(\tilde{\Omega}_{\varepsilon_k})} \le o(1)$ for $\varepsilon_k \to 0$, which is an absurd.

Up to now we have proved that:
\begin{equation}\label{HatPhiBoundLs}
\| \hat{\phi}_{\varepsilon} \|_{L^{2}(\tilde{\Omega}_{\varepsilon})} \le C
\end{equation}
To obtain an analogous bound on $\| \hat{\phi}_{\varepsilon} \|_{H^2(\tilde{\Omega}_{\varepsilon})}$ we can use  \cite[Lemma 6.4]{NiWei:Spike95} on $\hat{\phi}_{\varepsilon}$ together with Lemma \ref{PhiHatRHS} and \eqref{HatPhiBoundLs}. This concludes the proof. 
\end{proof}
Now the proof of Proposition \ref{PhiErrorEstimate} follows easily:
\begin{proof}
Since we have $\|\tilde u_{\varepsilon} \|_{H^2(\tilde{\Omega}_{\varepsilon})} \le C$, combining this, Lemma \ref{Lemma:PhiHatW2Estimate} and \eqref{PhiEstimateWRTPhiTilde} we get the assertion.
\end{proof}
Now we proceed with the estimate from below.
\begin{proposition}\label{ErrorEstimate}
For any sequence $\varepsilon_n \to 0$ there is a subsequence $\varepsilon_{n_k} \to 0$ such that:
\begin{equation}\label{EstimateFromBelow}
\tilde{\lambda}_{\varepsilon_{n_k}} = \tilde{\lambda}_0 + \Phi \, e^{-\beta_{\varepsilon_{n_k}} \tilde{\Psi}_{\varepsilon_{n_k}}(\mathbf{x}_{\varepsilon_{n_k}})} + o\biggl(e^{-\beta_{\varepsilon_{n_k}} \tilde{\Psi}_{\varepsilon_{n_k}}(\mathbf{x}_{\varepsilon_{n_k}})} \biggr) \quad \text{for } k \to +\infty \; ,
\end{equation}
where 
\[
\Phi \coloneqq \frac{2 \gamma}{\int_{\R^n} \tilde{m}_0 \, w^2}
\]
and $\gamma \coloneqq \int_{\R^n} f_0(w) \, \tilde{V}_0 > 0$ is independent on the particular $\tilde{V}_0$ solution of \eqref{eq_systemV0}
\end{proposition}
\begin{proof}
Using relation \eqref{eq:phi} combined with Lemma \ref{Lemma:PhiHatW2Estimate}, we have
\begin{align*}
& \tilde{\lambda}_{\varepsilon} = \frac{\int_{\tilde{\Omega}_{\varepsilon}} | \nabla \tilde u_{\varepsilon}|^2}{\int_{\tilde{\Omega}_{\varepsilon}} \tilde{m}_0 \, |\tilde u_{\varepsilon}|^2 } = \\
& = \frac{\int_{\tilde{\Omega}_{\varepsilon}} | \nabla P_{\tilde{\Omega}_{\varepsilon}} w|^2 + 2 e^{-\beta_{\varepsilon} \tilde{\Psi}_{\varepsilon}(\mathbf{x}_{\varepsilon})} \nabla P_{\tilde{\Omega}_{\varepsilon}} w \nabla \phi_{\varepsilon} + e^{-2\beta_{\varepsilon} \tilde{\Psi}_{\varepsilon}(\mathbf{x}_{\varepsilon})} |\nabla \phi_{\varepsilon}|^2
}{\int_{\tilde{\Omega}_{\varepsilon}} \tilde{m}_0 \, ((P_{\tilde{\Omega}_{\varepsilon}}w)^2 + 2 e^{-\beta_{\varepsilon} \tilde{\Psi}_{\varepsilon}(\mathbf{x}_{\varepsilon})} P_{\tilde{\Omega}_{\varepsilon}}w \phi_{\varepsilon} + e^{-2\beta_{\varepsilon} \tilde{\Psi}_{\varepsilon}(\mathbf{x}_{\varepsilon})}\phi_{\varepsilon}^2) } = \\
& = \frac{\int_{\tilde{\Omega}_{\varepsilon}} | \nabla P_{\tilde{\Omega}_{\varepsilon}} w|^2 + 2 e^{-\beta_{\varepsilon} \tilde{\Psi}_{\varepsilon}(\mathbf{x}_{\varepsilon})} \nabla P_{\tilde{\Omega}_{\varepsilon}} w \nabla \phi_{\varepsilon} + o(e^{-\beta_{\varepsilon} \tilde{\Psi}_{\varepsilon}(\mathbf{x}_{\varepsilon})})
}{\int_{\tilde{\Omega}_{\varepsilon}} \tilde{m}_0 \, ((P_{\tilde{\Omega}_{\varepsilon}}w)^2 + 2 e^{-\beta_{\varepsilon} \tilde{\Psi}_{\varepsilon}(\mathbf{x}_{\varepsilon})} P_{\tilde{\Omega}_{\varepsilon}}w \phi_{\varepsilon} + o(e^{-\beta_{\varepsilon} \tilde{\Psi}_{\varepsilon}(\mathbf{x}_{\varepsilon})}))} = \\
& = \frac{\int_{\tilde{\Omega}_{\varepsilon}} | \nabla P_{\tilde{\Omega}_{\varepsilon}} w|^2}{\int_{\tilde{\Omega}_{\varepsilon}} \tilde{m}_0 \, (P_{\tilde{\Omega}_{\varepsilon}}w)^2} + \frac{2e^{-\beta_{\varepsilon} \tilde{\Psi}_{\varepsilon}(\mathbf{x}_{\varepsilon})}}{\int_{\tilde{\Omega}_{\varepsilon}} \tilde{m}_0 \, (P_{\tilde{\Omega}_{\varepsilon}}w)^2} \biggl[ \int_{\tilde{\Omega}_{\varepsilon}} \nabla P_{\tilde{\Omega}_{\varepsilon}} w \nabla \phi_{\varepsilon} - \tilde{\lambda}_0 \tilde{m}_0 P_{\tilde{\Omega}_{\varepsilon}} w  \phi_{\varepsilon} \biggr]   + o(e^{-\beta_{\varepsilon} \tilde{\Psi}_{\varepsilon}(\mathbf{x}_{\varepsilon})}) \\
&= \frac{\int_{\tilde{\Omega}_{\varepsilon}} | \nabla P_{\tilde{\Omega}_{\varepsilon}} w|^2}{\int_{\tilde{\Omega}_{\varepsilon}} \tilde{m}_0 \, (P_{\tilde{\Omega}_{\varepsilon}}w)^2} + o(e^{-\beta_{\varepsilon} \tilde{\Psi}_{\varepsilon}(\mathbf{x}_{\varepsilon})}) 
\end{align*}
where in the last passage we have used the fact that both $P_{\tilde{\Omega}_{\varepsilon}} w$ and $\phi_{\varepsilon}$ belong to $H_0^1(\tilde{\Omega}_{\varepsilon})$ and $P_{\tilde{\Omega}_{\varepsilon}} w$ solves \eqref{eq:Ptildew}, so that
\[
\int_{\tilde{\Omega}_{\varepsilon}} \nabla P_{\tilde{\Omega}_{\varepsilon}} w \nabla \phi_{\varepsilon} - \tilde{\lambda}_0 \tilde{m}_0 P_{\tilde{\Omega}_{\varepsilon}} w  \phi_{\varepsilon} = 
 o(e^{-\beta_{\varepsilon} \tilde{\Psi}_{\varepsilon}(\mathbf{x}_{\varepsilon})})  \int_{\tilde{\Omega}_{\varepsilon}} f_0(\tilde{V}_{\varepsilon}) \phi_{\varepsilon} = o(e^{-\beta_{\varepsilon} \tilde{\Psi}_{\varepsilon}(\mathbf{x}_{\varepsilon})}) 
\]
thanks to Lemmas \ref{LemmaVTildeepsilon} and \ref{Lemma:PhiHatW2Estimate}. Using 
Remark \ref{rem:proofabove} the proof is concluded.
\end{proof}
We are ready to conclude the proof of our main result.
\begin{proof}[End of the proof of Theorem \ref{theorem:TheoremAsymptExpansion}]
Taking into account Propositions \ref{PhiErrorEstimate} and \ref{ErrorEstimate}, 
we are left to prove that for any sequence $\varepsilon_n \to 0$ there is a subsequence $\varepsilon_{n_k} \to 0$ such that:
\begin{equation}\label{positioning}
d(\mathbf{x}_{\varepsilon_{n_k}}, \partial\Omega) \to d(\mathbf{q}, \partial\Omega) \qquad \text{for } k \to +\infty
\end{equation}
\begin{equation}\label{PsiTilde}
\tilde{\Psi}_{\varepsilon_{n_k}}(\mathbf{x}_{\varepsilon_{n_k}}) \to 2 \sqrt{\tilde{\lambda}_0 \underline{m}} \, d(\mathbf{q}, \partial\Omega) \qquad \text{for } \varepsilon_{k_l} \to 0
\end{equation}
where $\mathbf{q}$ is any point at maximum distance from the boundary.

Combining Lemma \ref{LemmaVTildeepsilon} (ii) with \eqref{EstimateFromAbove} and \eqref{EstimateFromBelow} we have that, for any $\sigma_0 > 0$ there is an $\varepsilon_0 > 0$ such that, for $\varepsilon \le \varepsilon_0$
\begin{equation}\label{eqn:boundPsiTilde}
\Psi_{\varepsilon}(\mathbf{q}) \le \tilde{\Psi}_{\varepsilon}(\mathbf{x}_{\varepsilon}) \le (2+\sigma_0) \sqrt{\tilde{\lambda}_0 \underline{m}} \, d(\mathbf{x}_{\varepsilon}, \partial\Omega)
\end{equation}
Moreover, by Lemma \ref{LemmaVepsilon} (i) the above chain of inequalities can be rewritten as
\[
2 \sqrt{\tilde{\lambda}_0 \underline{m}} \, d(\mathbf{q}, \partial\Omega) + o(1) \le \tilde{\Psi}_{\varepsilon}(\mathbf{x}_{\varepsilon}) \le (2+\sigma_0) \sqrt{\tilde{\lambda}_0 \underline{m}} \, d(\mathbf{x}_{\varepsilon}, \partial\Omega) \qquad \text{for } \varepsilon \to 0
\]
Of course, it also holds that
\[
d(\mathbf{x}_{\varepsilon}, \partial\Omega) \le d(\mathbf{q}, \partial\Omega)
\]
so that
\[
\frac{2}{2+\sigma_0} d(\mathbf{q}, \partial\Omega) + o(1) \le  d(\mathbf{x}_{\varepsilon}, \partial\Omega) \le d(\mathbf{q}, \partial\Omega).
\]
This implies
\[
d(\mathbf{x}_{\varepsilon}, \partial\Omega) \to d(\mathbf{q}, \partial\Omega) \qquad \text{for } \varepsilon \to 0
\]
Thus, from \eqref{eqn:boundPsiTilde} we also have
\[
\tilde{\Psi}_{\varepsilon}(\mathbf{x}_{\varepsilon}) \to 2 \sqrt{\tilde{\lambda}_0 \underline{m}} \, d(\mathbf{q}, \partial\Omega) \qquad \text{for } \varepsilon \to 0
\]
This concludes the proof.
\end{proof}

\section{Appendix}
In this appendix we recall, for the sake of completeness, some of the main ideas leading to Proposition \ref{result:ResultProblemOverRn}. In particular, we show that all non-negative functions $u \in H^1(\R^N)$ that attain $\tilde{\lambda}_0$ in \eqref{lambda0defRn} for some $m \in \mathcal{M}'$ possibly dependent on $u$, need to be radially symmetric with respect to some point $\mathbf{x} \in \R^N$. This fact is the key point in the proof of Proposition \ref{prop:PropAwayFromBoundary}.

We begin with the following lemmas:
\begin{lemma}\label{lemma:LemmaSuperlevelOptimal}
Let $u \in H^1(\R^N)$. The quantity 
\begin{equation}\label{eqn:defMuOfu}
\mu(u) \coloneqq \sup_{m \in \mathcal{M}'} \int_{\R^N}m u^2
\end{equation}
is actually a maximum and is attained for 
\[
m' \coloneqq \overline{m} \ind{E_u}-\underline{m}\ind{ \R^N \setminus E_u } \; ,
\]
with
\[
E_u \coloneqq \{ \mathbf{x} \in \R^N : u^2 \ge \alpha \}
\]
and 
\[
\alpha \coloneqq \inf_{a \ge 0} \{ a : \mathcal{L}(\{ u^2 \ge a  \}) < \mathcal{L}(B_1(\mathbf{0}) ) \} \; .
\]
The above definition of $E_u$ is valid if $\mathcal{L}(\{ u^2 = \alpha \}) = 0$. Otherwise, a measurable portion of the set $\{ \mathbf{x} \in \R^N : u^2 = \alpha \}$ has to be removed from $E_u$.
\end{lemma}
\begin{proof}
We follow the ideas in \cite{CantrellCosner:indefiniteWeight} Theorem 3.9.

Consider
\[
\alpha = \inf_{a \ge 0} \{ a : \mathcal{L}(\{ u^2 \ge a  \}) < \mathcal{L}(B_1(\mathbf{0}) ) \} \; .
\]
If $\mathcal{L}(\{ u^2 = \alpha \}) = 0$, then $\alpha > 0$ and $E_u = \{ \mathbf{x} \in \R^N : u^2 \ge \alpha \}$. 

Otherwise, consider eventually just a measurable subset $\mathcal{I}$ of $\{\mathbf{x} \in \R^N : u^2 = \alpha \}$ chosen such that $\mathcal{L}(\{  u^2 > \alpha \} \cup \mathcal{I}) = \mathcal{L}(B_1(\mathbf{0}) )$. In this case we set $E_u = \{  u^2 > \alpha \} \cup \mathcal{I}$. We remark that in this second case we may also have $\alpha = 0$. 

We will only give the proof for the first case, since the second one is similar. Define
\[
m'(\mathbf{x}) \coloneqq
\begin{cases}
\overline{m} & \text{if } \mathbf{x} \in E_u \; ,\\
\underline{m} & \text{if } \mathbf{x} \in \R^N \setminus E_u \; ,
\end{cases}
\]
and let $m \in \mathcal{M}'$. Then:
\begin{align*}
& \int_{\R^N} (m' - m) u^2 = \int_{\R^N} [(m'+\underline{m}) - (m+\underline{m})] u^2 = \\
& = \int_{E_u} [(m'+\underline{m}) - (m+\underline{m})] u^2 + \int_{\R^N \setminus E_u} [(m'+\underline{m}) - (m+\underline{m})] u^2 \ge \\
& \ge \alpha \int_{E_u} [(m'+\underline{m}) - (m+\underline{m})] - \alpha \int_{\R^N \setminus E_u} [(m+\underline{m}) - (m'+\underline{m})] = \\
& = \alpha \int_{\R^N} [(m'+\underline{m}) - (m+\underline{m})] = \alpha \left[ (m'+\underline{m}) \mathcal{L}(B_1(\mathbf{0}) ) - \int_{\R^N}(m+\underline{m}) \right] \ge 0
\end{align*}
where in the last passage we have used the fact that, by the definition of $\mathcal{M}'$, we have
\[
\int_{\R^N}(m+\underline{m}) \le (\underline{m} + \overline{m}) \mathcal{L}(B_1(\mathbf{0}) ) \; .
\]
Thus we have proved that
\[
\int_{\R^N} (m' - m) u^2 \ge 0
\]
for any $m \in \mathcal{M}'$ and this concludes the proof.
\end{proof}
\begin{lemma}\label{lemma:LemmaSuperlevelSetMeasure}
Let $\mathcal{A} \subseteq \R^N$ be an open set and $f \in C^2(\mathcal{A})$. Define the set
\[
\mathcal{C} \coloneqq \{ \mathbf{x} \in \mathcal{A} : \nabla f = 0 \} \; .
\]
Then, for almost every $\mathbf{x} \in \mathcal{C}$ it holds that $\partial_{x_i x_j} f = 0$, for any $i, j = 1, ..., N$.
\end{lemma}
\begin{proof}
It is a standard fact in measure theory (for instance \cite{Tao:introMeasuteTheory}, p. 153) that if $E \in \R^N$ is a Lebesgue measurable set, then almost every point of $E$ is a point of density. This means that:
\begin{equation}\label{eqn:densityPropertyOfMeasurableSet}
\frac{\mathcal{L}(E \cap B_r(\mathbf{x}))}{\mathcal{L}(B_r(\mathbf{x}))} \to 1 \qquad \text{for } r \to 0^{+}
\end{equation}
for almost every $\mathbf{x} \in E$.

Consider $E = \mathcal{C}$. If $\mathcal{L}(\mathcal{C}) = 0$ then we have nothing more to prove. Thus, suppose $\mathcal{L}(\mathcal{C}) > 0$. We will prove that in each point of density for $\mathcal{C}$ it holds that $\partial_{x_i x_j} f = 0$.

Fix a point $\mathbf{x} \in \mathcal{C}$ of density for $\mathcal{C}$. For any direction $i$ consider a half-cone $K_{\mathbf{x}, \Omega}$ centered at $\mathbf{x}$ of arbitrary small solid angle $\Omega$, whose axis has direction $i$. Thanks to \eqref{eqn:densityPropertyOfMeasurableSet}, for any $r > 0$ sufficiently small, we have that $\mathcal{I}_r \coloneqq \mathcal{C} \cap B_r(\mathbf{x}) \cap K_{\mathbf{x}, \Omega} \neq \emptyset$. Consider a sequence $r_n \to 0$ for $n \to +\infty$ and $\mathbf{x}'_n \in \mathcal{I}_{r_n}$. Then, by definition, for any $\mathbf{x}'_n$ we have $\partial_{x_j} f(\mathbf{x}'_n)$ = 0, for any $j$. Hence, Lagrange theorem implies that, for $n$ large enough so that $B_{r_n}(\mathbf{x}) \subset \mathcal{A}$, $\frac{\partial}{\partial\nu_n} \partial_{x_j} f(\mathbf{x}_n^{\ast}) = 0$. Here  $\mathbf{x}_n^{\ast}$ is a point lying on the segment connecting $\mathbf{x}$ to $\mathbf{x}'_n$, while $\nu_n$ denotes its direction, which lies inside $K_{\mathbf{x}, \Omega}$.

Now, since $\frac{\partial}{\partial\nu_n} \partial_{x_j} f(\mathbf{x}_n^{\ast})$ can be written as $\sum_k \partial_{x_k x_j} f(\mathbf{x}_n^{\ast}) \nu_{n, k}$, where $\nu_{n, k}$ denotes the k-th component of the normal $\nu_n$, considering that all the second derivatives are bounded since $f \in C^2(\bar{B}_{r_n}(\mathbf{x}))$, we have
\[
\partial_{x_i x_j} f(\mathbf{x}) = o(1) \qquad \text{for } \Omega \to 0^{+}, \; n \to +\infty\; ,
\]
because $\mathbf{x}_n^{\ast} \to \mathbf{x}$ for $n \to +\infty$ and the second partial derivatives are continuous, moreover $\nu_{n, k} \to 0$ if $k \neq i$ for $\Omega \to 0^{+}$ uniformly in $n$ and $\nu_{n, i} \to \pm 1$ for $\Omega \to 0^{+}$ uniformly in $n$, depending on the orientation of the half-cones. 

Thus, letting $\Omega \to 0^{+}$ and $n \to +\infty$ we get
\[
\partial_{x_i x_j} f(\mathbf{x}) = 0
\]
for any point of concentration of $\mathcal{C}$, and so for almost every point in $\mathcal{C}$. This concludes the proof.
\end{proof}
Now consider a non-negative function $u \in H^1(\R^N)$ attain $\tilde{\lambda}_0$ in \eqref{lambda0defRn}. We are ready to show that it needs to be radially symmetric with respect to some $\mathbf{x} \in \R^N$.

Indeed, consider $u$ normalized such that $\int_{\R^N}|\nabla u|^2 = 1$. Then, we can write
\[
\tilde{\lambda}_0 = \frac{1}{\int_{\R^N}m u^2} \; .
\]
However, thanks to Lemma \ref{lemma:LemmaSuperlevelOptimal} we can find $m' \in \mathcal{M}'$ for which $\mu(u)$ is attained in \eqref{eqn:defMuOfu}. Thanks to the definition of $m'$ it holds that
\[
\frac{1}{\int_{\R^N}m' u^2} \le  \frac{1}{\int_{\R^N}m u^2} \; .
\]
But thanks to the normalization chosen for $u$ and to the minimality property \eqref{lambda0defRn} of $\tilde{\lambda}_0$ we have
\[
\tilde{\lambda}_0 \le \frac{1}{\int_{\R^N}m' u^2}  \le \frac{1}{\int_{\R^N}m u^2} = \tilde{\lambda}_0 \; ,
\]
and thus
\[
\frac{1}{\int_{\R^N}m' u^2} = \tilde{\lambda}_0 \; .
\]
Now consider the normalization for $u$ such that $\int_{\R^N} m' u^2 = 1$. We still denote it as $u$. If we introduce the Schwarz symmetrization $u_{\ast}$ of $ u $, thanks to the superlevel properties of $m'$ remarked in Lemma \ref{lemma:LemmaSuperlevelOptimal}, we have:
\[
1 = \int_{\R^N} m' u^2 = \int_{\R^N} m'_{\ast} u_{\ast}^2 \; ,
\]
where $m'_{\ast} = \tilde{m}_0$ is the spherical symmetrization of $m'$.

By the P\'olya-Szeg\H{o} inequality in \cite{PolyaSzego:isoperimetricInequalitites,Henrot:extremumProblems} then we have:
\[
\int_{\R^N}|\nabla u_{\ast}|^2 \le \int_{\R^N}|\nabla u|^2 = \tilde{\lambda}_0 \; ,
\]
and by the minimality property \eqref{lambda0defRn} satisfied by $\tilde{\lambda}_0$, it holds that
\begin{equation}\label{eqn:PolyaSzegoEqualityContradictionDistanceBoundary}
\int_{\R^N}|\nabla u_{\ast}|^2 = \int_{\R^N}|\nabla u|^2 \; .
\end{equation}
We can conclude by the equality case in \cite[Thm. 1.1]{BrothersZiemer:minimalRearrangements}, whose hypotheses are satisfied. Indeed, thanks to \eqref{eqn:PolyaSzegoEqualityContradictionDistanceBoundary}, taking Gateaux derivatives of the functional
\[
\frac{\int_{\R^N}|\nabla u_{\ast}|^2}{\int_{\R^N}\tilde{m}_0 u_{\ast}^2}
\]
and exploiting the minimality of $\tilde{\lambda}_0$, it follows that $u_{\ast}$ solves, in $H^1(\R^N)$, the equation
\[
-\Delta u_{\ast} = \tilde{\lambda}_0 \tilde{m}_0 u_{\ast} \quad \text{in } \R^N \, .
\]
Thanks to Lemma \ref{lemma:LemmaSuperlevelSetMeasure}, this implies that
\begin{equation}\label{eqn:setMeasureEqualityContradictionDistanceBoundary}
\mathcal{L}(\mathcal{A}) = 0, \text{ where }\mathcal{A} \coloneqq \{ \mathbf{x} \in \R^N : u_{\ast}>0,  \nabla u_{\ast} = 0\} \; .
\end{equation}
Indeed, suppose not. By standard elliptic regularity we have that in $B^1(\mathbf{0})$ and $\R^N \setminus B^1(\mathbf{0})$ the function $u_{\ast}$ is $C^2$, and thus it holds that $\Delta u_{\ast} < 0$ in $\mathcal{A} \cap B^1(\mathbf{0})$ and $\Delta u_{\ast} > 0$ in $\mathcal{A} \cap B^1(\mathbf{0})^c$. This means that $\Delta u_{\ast} \neq 0$ almost everywhere in $\mathcal{A}$, and we have reached an absurd.

Thus, thanks to \cite[Thm. 1.1]{BrothersZiemer:minimalRearrangements},  \eqref{eqn:PolyaSzegoEqualityContradictionDistanceBoundary} and \eqref{eqn:setMeasureEqualityContradictionDistanceBoundary} we have that $u$ needs to be radially symmetric and actually also radially non increasing with respect to a point $\mathbf{x} \in \R^N$.

\section*{Acknowledgments} Work partially supported by the project Vain-Hopes within the program VALERE: VAn- viteLli pEr la RicErca, by the Portuguese government through FCT/Portugal under the project PTDC/MAT-PUR/1788/2020, and by the INdAM-GNAMPA group.

\bibliography{Ferreri_Bibliography}
\bibliographystyle{abbrv}

\medskip
\small
\begin{flushright}
\noindent 
\verb"lorenzo1.ferreri@mail.polimi.it"\\
\verb"gianmaria.verzini@polimi.it"\\
Dipartimento di Matematica, Politecnico di Milano\\ 
piazza Leonardo da Vinci 32, 20133 Milano (Italy)
\end{flushright}

\end{document}